\documentclass[11pt]{amsart}
 
\usepackage[cm]{fullpage}
\usepackage{appendix}
\usepackage{latexsym}
\usepackage{amsfonts}
\usepackage{amssymb}
\usepackage{tikz-cd}
\usepackage{enumitem}
\usepackage{hyperref}
\usepackage{graphicx} 
\graphicspath{ {./images/} }
\usepackage{comment}

\usepackage{amsmath,amssymb,color}
\usepackage{amsthm}
\usepackage{hyperref}
\usepackage{xcolor}

\usepackage{amsmath,amscd, float,rotating}
\usepackage{pb-diagram}

\numberwithin{equation}{section}

\def\XXint#1#2#3{{\setbox0=\hbox{$#1{#2#3}{\int}$}
\vcenter{\hbox{$#2#3$}}\kern-.5\wd0}}

\def\na{\nabla}

\newcommand{\al}{\alpha}
\newcommand{\be}{\beta}

\newcommand{\vep}{\varepsilon}

\newcommand{\cU}{\mathcal{U}}

\newcommand{\KE}{\mathrm{KE}}

\newcommand{\dbar}{\overline{\partial}}

\newcommand{\ddbar}{\sqrt{-1}\partial\dbar}

\newtheorem{theorem}{Theorem}[section]
\newtheorem{proposition}{Proposition}[section]
\newtheorem{lemma}{Lemma}[section]
\newtheorem{example}{Example}[section]

\newtheorem{conjecture}{Conjecture}[section]
\newtheorem{definition}{Definition}[section]
\newtheorem{corollary}{Corollary}[section]
\newtheorem{remark}{Remark}[section]

\begin{document}

~

\title{K\"ahler-Einstein metrics  near an isolated log-canonical singularity}

\author{Ved Datar$^*$, Xin Fu$^{**}$, Jian Song$^\dagger$ }

\address{$*$ Department of Mathematics, Indian Institute of Science, Bangalore, India}

\email{vvdatar@iisc.ac.in}

\address{$**$ Department of Mathematics,  University of California, Irvine, CA 92617}

\email{fux6@uci.edu}

\address{$\dagger$ Department of Mathematics, Rutgers University, Piscataway, NJ 08854}

\email{jiansong@math.rutgers.edu}

\thanks{Research supported in part by National Science Foundation grant  DMS-1711439, and by the Infosys Young Investigator award}

\begin{abstract} We construct  K\"ahler-Einstein metrics with negative scalar curvature near an isolated log canonical (non-log terminal) singularity. Such metrics are complete near the singularity if the underlying space has complex dimension $2$. We also establish a stability result for K\"ahler-Einstein metrics near certain types of isolated log canonical singularity. As application, for complex dimension $2$ log canonical singularity, we show that any complete K\"ahler-Einstein metric of negative scalar curvature near an isolated log canonical (non-log terminal) singularity is smoothly asymptotically close to model K\"ahler-Einstein metrics from hyperbolic geometry.

\end{abstract}

\maketitle

%%%%%%%%%%%%%%%%%%%%%%%%%%%%%%%%%%%%%%%%

\section{Introduction}\label{intro}

The existence  of K\"ahler-Einstein metrics on complex manifold has been the central topic in complex geometry for decades. In \cite{Y1}, Yau  established the existence of Ricci flat metrics on complex manifolds with zero Chern class by solving the Calabi conjecture, while Aubin and Yau \cite{A,Y1} proved the existence of K\"ahler-Einstein metrics independently on canonically polarized compact manifolds.  
Recent results of Chen-Donaldson-Sun \cite{CDS1,CDS2,CDS3}, Tian \cite{T4}  and many interesting subsequent papers proved equivalence between  existence of K\"ahler-Einstein metrics and $K$-stability for Fano manifolds, confirming the  Yau-Tian-Donaldson conjecture. 
Also there has been intensive study of degenerate complex Monge-Amp\`ere equations and construction of singular K\"ahler-Einstein metrics on singular varieties with Kawamata log terminal (klt) singularities, for example in \cite{EGZ,Z},  based on Kolodziej's fundamental result in \cite{Kol1}. K\"ahler-Einstein metrics on canonical polarized variety with log canonical singularity are constructed in  \cite{B} by the variational approach. When such varieties appear on the boundary of the KSB compactification of smooth canonical models of general type, such K\"ahler-Einstein metrics turn out to be genuinely geometric as degeneration of K\"ahler-Einstein metrics \cite{S4, SSW}.  However, little is known about the geometric behavior of these singular K\"ahler-Einstein metrics near the log canonical singularities. For smoothable klt singularities, the fundamental work of Donaldson-Sun \cite{DS1,DS2} and \cite{LX} show that the tangent cone at any given singular point  is  unique and admits a Ricci-flat cone metric as the Gromov-Hausdorff limit of smooth K\"aher-Einstein metrics.  Such non-collapsed local models for smoothable klt singularities are always built on Sasakian geometry and an interesting metric stability result is obtained   in \cite{HS,CS}. In this paper, we will study  K\"ahler-Einstein metrics near isolated log canonical (non-log terminal) singularity and their asymptotic collapsing behavior near the complete end. 

%We conjecture that any complete K\"ahler-Einstein metrics near log canonical singularities are asymptotically close to each other with a local collapsing model of suitable scalings.     

Let $(X,p)$ be a germ of an isolated log canonical algebraic singularity $p$ embedded in $(\mathbb C^N,0)$. In this paper, we will assume $p$ is not a log terminal singularity. We would like to  construct local K\"ahler-Einstein metrics with negative scalar curvature in an open neighborhood of the singular point $p$. 
%Since the local volume of $p$ by generators of the canonical divisor is not integrable, it seems more appropriate to discuss K\"ahler-Einstein metric with negative scalar curvature near the singular point $p$. 
Let $\rho$ be a  nonnegative smooth plurisubharmonic (PSH) function on $\mathbb C^N$ with $\rho(0)=0$ and $\ddbar \rho>0$.  We let  $\cU \subset\subset X$ be the open domain defined by
\begin{equation}\label{domain1}
\cU=\{\rho< a\} \cap X
\end{equation}
for some sufficiently small constant $a>0$. 
For example, we can always choose $\rho = \sum_{j=1}^N |z_j|^2$. By choosing a  generic sufficiently small $a>0$, we can assume that $\partial\cU$ is smooth and strongly pseudoconvex. We will fix such a domain $\cU$ for $(X, p)$ for the rest of the paper. We further require that the canonical divisor $K_{\cU}$ is Cartier and fix a local volume measure $\Omega_X$ on $\cU$ by 
\begin{equation}\label{volumeform}\Omega_X=  (\sqrt{-1})^n\nu \wedge \overline \nu, \end{equation} 
where $\nu$ is a local generator of $K_\cU$. Similarly, we can define $\Omega_X$ when $K_{\cU}$ is $\mathbb{Q}$-Cartier. 

We will consider the following Dirichlet problem of a complex Monge-Amp\`ere equation related to the K\"ahler-Einstein metric on $\cU$.
\begin{equation} \label{ocmp}
\left\{
\begin{array}{ll}
(\ddbar \varphi)^n = e^\varphi \Omega_X, &  in ~ \cU\setminus\{p\} \\
& \\ 
\varphi =\psi, & on ~ \partial \cU.\\
\end{array}\right. 
\end{equation}
where $\psi\in C^\infty(\partial \cU)$. Our first result is on existence of finite volume K\"ahler-Einstein metrics near isolated log-canonical singularities. To state the theorem, fix a log resolution $\pi:Y\to X$ and let $D$ be the simple normal crossing exceptional divisor. We also fix a defining section $\sigma_D$ and a hermitian metric $h_D$ for the line bundle corresponding to $D$ (cf. Section 2 for more precise definition). We 
\begin{theorem}\label{main1} 
Let $(\cU,p)$ be a germ of an isolated  log canonical algebraic (non-log terminal) singularity defined as above. For any smooth function $\psi$ on $\partial\cU$, there exists a function $\varphi_{\KE}$  satisfying the following:
\begin{enumerate}

\item $\varphi_{\KE}\in C^{\infty}(\pi^{-1}(\cU)\setminus \mathrm{Supp}(D)) \cap PSH(\pi^{-1}(\cU))$, 

 \item For any $\epsilon>0$, there is a constant $C_\epsilon>0$, such that on $\pi^{-1}(\cU)\setminus \mathrm{Supp}(D)$,$$C\geq \varphi_{\KE} \geq -(2n+\epsilon)
 \log(-\log |\sigma_D|^2_{h_D}) - C_\epsilon.$$
\item  $\varphi_{\KE}$ solves equation \eqref{ocmp} on $\cU\setminus \{p\}$ if we identify $\cU\setminus\{p\}$ with $\pi^{-1}(\cU)\setminus\mathrm{Supp}(D)$, and the K\"ahler metric $\omega_\KE = \ddbar\varphi_\KE$ is K\"ahler-Einstein with finite volume: $$\int_{\cU\setminus\{p\}}\omega_\KE^n<\infty.$$
\item $\varphi_{\KE}$ has vanishing Lelong number at any point in $\pi^{-1}(\cU)$ and for $x\in\pi^{-1}(\cU)$,$$\lim_{x\to \mathrm{Supp}(D)}\varphi_{\KE}(x)= -\infty.$$

\end{enumerate}
Moreover, any solution to \eqref{ocmp} (when 
 pulling back to $\pi^{-1}(\cU)\setminus\mathrm{Supp}(D)$) satisfying $(1)$ and $(2)$ above is unique.

\end{theorem}

%\begin{remark}
%In item(4) above, $\varphi_{KE}$ has vanishing Lelong number means that $\Pi^*\varphi_{KE}$ has vanishing Lelong number as a PSH function on smooth domain $Y$. 
%\end{remark}

%Since $\varphi$ has vanishing Lelong number,  $g_{KE}$ can be uniquely extended to K\"ahler currents on $\cU$. Therefore equation \eqref{keq} also holds on $\cU$ as currents. In fact, it can shown by similar estimates in \cite{} (c.f. \cite[Proposition 3.1]{SSW}) that $\varphi$ can have at worst $\log(-\log)$-poles, and so $Ric(g_{KE})$ can also be extended on $\cU$ as a negative current. 

One naturally would ask how the K\"ahler-Einstein metric $g_{KE}$ in Theorem \ref{main1} behaves near the singularity $p$ geometrically. We propose the following conjecture.

\begin{conjecture} \label{conj1} Let $g_{KE}$ be a K\"ahler-Einstein metric constructed in Theorem \ref{main1} on $\cU \setminus \{p\}$. Then $g_{KE}$ is complete  near $p$. More precisely, for any fixed point $q\in \cU\setminus\{p\}$ and $R>0$, there exists an open neighborhood $V_{q, R}$ of $p$ in $\cU$ such that for any point $x \in V_{q, R} \setminus \{p\}$ and any smooth path $\gamma$ joining $q$ and $x$ in $\cU \setminus \{p\}$, we have 
$$|\gamma|_{g_{KE}} >R,$$
where $|\gamma|_{g_{KE}}$ is the arc length of $\gamma$ with respect to $g_{KE}$.

\end{conjecture}

We are able to confirm the above conjecture in the following two special cases. 

\begin{theorem}\label{conj} Conjecture \ref{conj1} holds when $\dim \cU=2$.  

\end{theorem}

The construction from Theorem \ref{main1} gives infinitely many K\"ahler-Einstein metrics near the log canonical singularity $p$ by assigning different boundary conditions.  However, certain asymptotic stability should hold for such complete K\"ahler-Einstein metrics.  More precisely,  the asymptotic geometric behavior of complete K\"ahler-Einstein metrics near a log canonical singularity $p$ should be completely determined by the analytic structure of $p$.

%\label{setting 2} Let ($X,p$) be germ of isolated log canonical singularity, if there is K\"ahler-Einstein metric $\omega_{KE}=\ddbar\rho$  with bounded geometry property (see (\ref{property A}) for details) on a punctured neighbourhood of $p$ and complete  towards $p$, we define $$\Omega:=\{\rho<a\},\partial\cU:=\{\rho=a\}$$ where $a$ is a fixed constant. We remark here that the behaviour of $\rho$ here is quite different from the $\rho$ we used to define domain $\Omega$ in theorem (\ref{main1}). 
%\begin{definition}
% Complete metric: In this note, when we say a K\"ahler metric is complete, we mean that it is complete with infinite distance towards $p$ and incomplete with finite distance towards $\partial\cU$. 

%\end{definition}

Our next result  is to  establish a volume stability for such complete K\"ahler-Einstein metrics. We can always shrink  $\cU$ slightly because we are interested in the behavior of complete K\"ahler-Einstein metrics near the isolated log canonical singularity $p$.  

\begin{theorem}\label{volume} Let  $g_{KE}$ and $g'_{KE}$ be two K\"ahler-Einstein metrics on $\cU \setminus \{ p\}$ for an isolated log canonical (non-log terminal) singularity $p$.  Let $R(x)$ and $R'(x)$ be the distance functions from any point $x\in \cU\setminus\{p\}$ to $\partial \cU$ with respect to $g_{KE}$ and $g'_{KE}$.  
If  both $g_{KE}$ and $g'_{KE}$ are complete near $p$, then there exits $c=c(n)>0$ such that 
\begin{equation}1-\frac{c(n)}{R'(x)}\le
\frac{\det \left( g'_{KE} \right) }{ \det \left( g_{KE} \right) }(x)\le 1+\frac{c(n)}{R(x)}\end{equation}
for any $x \in \cU \setminus\{p\}$ with $R(x)>1$ and $R'(x)>1$. In particular, 
$$\lim_{x\rightarrow p} \frac{\det \left(g'_{KE}\right)}{\det \left( g_{KE} \right)}(x) = 1,$$
where $x \rightarrow p$ with respect to the Euclidean metric from any local embedding of $p$ in $\mathbb{C}^N$.
\end{theorem}
Theorem \ref{volume} also implies that the potentials of two complete K\"ahler-Einstein metrics near the isolated log canonical singularity $p$ must be asymptotically close to each other at infinity.

There is a large class of log canonical singularities that admits complete K\"ahler-Einstein metrics with bounded geometry at the complete end. For example, a metric uniformization  is obtained in \cite{Kob,Kob2} for  isolated log canonical singularity in complex dimension $2$. Another interesting example of cones over an abelian variety is constructed in \cite{Fang-Fu}. We apply the method of bounded geometry
by \cite{CY} (see also \cite{Kob,Kob2,TY}) to prove the following stability result.

\begin{theorem}\label{main2} Let $(\cU, p)$ be a germ of an isolated log canonical singularity $p$. If there exists a complete K\"ahler-Einstein metric $\theta$ near $p$ with bounded geometry of order $k\geq 2$, then for any complete K\"ahler-Einstein metric $g_{KE}$ on  $\cU\setminus \{p\}$ and $k \geq 0$, there exists $C= C(n,k,\theta, g_{KE})>0$ such that for any $x\in \cU$ with $R_\theta(x)>1$, we have
 
\begin{equation}
    \|\nabla_\theta^k(\omega_{KE} - \theta )\|_{\theta}(x)\le \frac{C}{\left( R_\theta(x) \right)^{1\over 2}}
    ,
    \end{equation}
 where $R_\theta(x)$ is the distance from $x$ to $\partial U$ and $\nabla_\theta$ is the covariant derivative with respect to $\theta$.

\end{theorem}
%\begin{remark}
%We wish to now compare our result to classical results obtained by using the bounded geometry method.  This construction of complete K\"ahler-Einstein metrics using the method of bounded geometry was initiated by the fundamental paper of Cheng-Yau \cite{CY}.  These classical results imply that the constructed K\"ahler-Einstein is equivalent to a model reference metric at infinity. Our improvements in Theorem \ref{main2} are twofolds :  Firstly, and most crucially, we show that any complete local K\"ahler-Einstein metrics are equivalent to each other (and hence all are equivalent to any  model metric that might exist).   To our knowledge, this was not known before. Secondly, we show that the complete local K\"ahler-Einstein is in fact {\em asymptotic} to the model metric. \end{remark}
\begin{remark}
We  remark that the completeness assumption of $\omega_{KE}$ is necessary, because there are examples of germ of log canonical singularity $(\cU\setminus p)$ which admit both complete and incomplete K\"ahler-Einstein metrics (cf. \cite[Example 2.7]{FHJ}).  We also want to compare our stability result with the results of \cite{HS,CS}. There in the Ricci flat case, if we could relate two different K\"ahler-Einstein metrics $\omega^1_{KE},\omega^2_{KE}$ as $\omega^1_{KE}=\omega^2_{KE}+\ddbar u$, the difficulty lies in the higher order estimate of $u$. While in our case, one new ingredient is that we can show the boundedness of volume ratio $u$ unconditionally in Theorem \ref{volume} and then use it to show that any local complete $\omega_{KE}$  can be obtained from Cheng-Yau type construction once we have a good model metric. Then higher order regularity of volume ratio is a byproduct of Cheng-Yau type construction. And to our knowledge, the asymptotics result improves Cheng-Yau global construction and hence is new. 
\end{remark}

We also have the following immediate corollary.

\begin{corollary} \label{surface}Let $(\cU, p)$ be the germ of a log canonical singularity $p$ with $\dim \cU=2$. Then for any complete K\"ahler-Einstein metric $g_{KE}$ near $p$, $(\cU\setminus\{p\}, g_{KE})$ must be asymptotically isometric to one of the following two local models.

\begin{enumerate}

\item  $\left( \mathbb{B}^2 /\Gamma, g_{\mathbb{B}^2 /\Gamma} \right)$, where $\mathbb B^2$ is the unit ball in $\mathbb{C}^2$, $\Gamma$ is a parabolic discrete subgroup of $\textnormal{Aut}(\mathbb{B}^2)$ and $g_{\mathbb{B}^2 / \Gamma} $ is the  K\"ahler metric induced by the group arction invariant hyperbolic metric on $\mathbb{B}^2$.

\item  $\left((\mathbb{H}\times \mathbb{H} )/  \Gamma, g_{(\mathbb{H} \times \mathbb{H}) /  \Gamma} \right)$, where $\mathbb{H}$ is the hyperbolic upper half plane, $\Gamma$ is a parabolic discrete subgroup of $\textnormal{Aut}(\mathbb{H}\times \mathbb{H})$ and $g_{(\mathbb{H} \times \mathbb{H}) /  \Gamma}$ is the  K\"ahler metric induced by the group action invariant hyperbolic metric $g_{\mathbb{H}\times\mathbb{H}}$.
\end{enumerate}
Furthermore, if we take a sequence of points $p_j\in (\cU\setminus p)$ with  $p_j \rightarrow p$ (in the Euclidean topology), then $(\cU, p_j, g_{KE})$ converges to line $\mathbb{R}$ as in case \textnormal{(1)} and converges to  flat cylinder $S^1\times \mathbb{R}$ as in case \textnormal{(2)},  in pointed Gromov-Hausdorff topology.

\end{corollary}

Theorem \ref{main2} can also be applied to higher dimensional complex hyperbolic cusps. (cf. \cite{Fang-Fu})

\begin{remark} For application, our Theorem \ref{main2} and Corollary \ref{surface} can be applied to understand the geometry of K\"ahler-Einstein metric with negative scalar curvature constructed on certain stable canonical polarized variety by Berman-Guenancia \cite{B} and Song \cite{S4} (cf. example \ref{ex}).
\end{remark}
We further conjecture that for any complete K\"ahler-Einstein metric $g_{KE}$ on $(\cU, p)$   an isolated log canonical singularity $p$, any sequence of points $x_j \rightarrow p$ and positive $\lambda_j \geq 1$, then $(\cU \setminus \{p\}, x_j, \lambda_j g_{KE})$ converges in pointed Gromov-Hausdorff topology to a product of $\mathbb{R}$, $\mathbb{C}$, compact Calabi-Yau varieties and complete Calabi-Yau varieties with cylindrical end.

We briefly outline our paper. In section 2, we prove Theorem \ref{main1} by solving the Dirichlet problem of singular complex Monge-Ampere equations. We prove the volume stability of Theorem \ref{volume} for complete K\"ahler-Einstein metrics near isolated log canonical singularities in section 3 and a metric stability  Theorem \ref{main2} if one of the K\"ahler-Einstein metric satisfies the bounded geometry condition in section 4. In section 5, we prove Theorem \ref{conj} for surfaces and Corollary \ref{surface}. In section 6, we give a short discussion of conjecture \ref{conj1} for smoothable isolated log singularities.
\bigskip

\noindent {\bf{Acknowledgements:}} The second author would like to thank Professor H-J. Hein for his illuminating lectures on his joint work with Professor S. Sun \cite{HS} during the Summer School at Notre Dame University in 2017 and  F. Tong for many interesting discussions. We also thank the anonymous referee for many useful suggestions.
 
\section{K\"ahler-Einstein metrics near log canonical singularities}

\subsection{The set-up.}\label{sec:set-up} In this section, we will prove Theorem \ref{main1}. We first recall the definition for log canonical singularities.

\begin{definition} \label{sing} Let $X$ be a normal  variety such that $K_X$ is a $\mathbb{Q}$-Cartier divisor. Let $ \pi :  Y \rightarrow X$ be a log resolution and $\{E_i\}_{i=1}^p$ the irreducible components of the exceptional locus $Exc(\pi)$ of $\pi$. There there exists a unique collection $a_i\in \mathbb{Q}$ such that 
$$K_Y = \pi^* K_X + \sum_{i=1}^{ p } a_i E_i .$$ Then $X$ is said to have log-canonical (resp.\ klt) singularities if $$a_i\geq -1~ \text{(resp.\ $a_i>-1$) for all $i$}. $$
%
%\begin{enumerate}

%\item[$\bullet$] terminal singularities if  $a_i >0$, for all $i$.

%\item[$\bullet$] canonical singularities if $a_i \geq 0$, for all $i $.

%\item[$\bullet$]  log terminal singularities if $a_i > -1$, for all $i $.

%\item[$\bullet$]  log canonical singularities if $a_i \geq -1$, for all $i$.

%\end{enumerate}
\end{definition}

In this paper, since we are considering isolated singularity, we have $\pi(Exc(\pi))={p}$.
We have to prescribe singularities of the solution $\varphi$ to obtain a canonical and unique K\"ahler-Einstein current on $X$. To do so, we   lift all the data to a log resolution $\pi: Y \rightarrow X$. 
It is more convenient to write the adjunction formula in the following form:
 $$K_Y = \pi^* K_X + \sum_i^n a_i E_i - \sum_j^m b_j F_j, ~a_i\geq 0, ~ 0< b_j \leq 1. $$

%We approximate  equation (\ref{ocmp})  in  the following way. We pull back all the data from $X$ to $Y$. 

Let $\sigma_{E_i}, 1\leq i\leq n$ be the defining section for line bundle associated to $E_i$ and $\sigma_{F_i}, 1\leq j\leq m$ be the defining section for line bundle associated to  $F_j$. We also equip the line bundles associated to $E_i$ and $F_i$ with smooth hermitian metric $h_{E_i}$, $h_{F_j}$ on $Y$. Then we define:
\begin{equation*}|\sigma_E|^2_{h_E}:=\prod_i|\sigma_{E_i}|_{h_i}^{2a_i},  |\sigma_F|^2_{h_F}:=\prod_j|\sigma_{F_j}|_{h_j}^{2b_j}\end{equation*}
Recall that $\Omega_X$ is a local volume form defined on $\cU$ in formula (\ref{volumeform}), now let $\Omega_Y$ be a smooth strictly positive volume form on $\pi^{-1}(\cU)$,  
defined by 
$$ \Omega_Y  = ( |\sigma_E|^2_{h_E} )^{-1} |\sigma_F|^2_{h_F} \pi^*\Omega_X. $$
Then lifting equation \eqref{ocmp} to $Y$, we have 

\begin{align}\label{aaaa}
\begin{cases}
&(\pi^*\ddbar \varphi)^n = e^{\pi^*\varphi} ( |\sigma_F|^2_{h_F} )^{-1} |\sigma_E|^2_{h_E} \Omega_Y, \\
&\pi^*\varphi_{|\partial\cU}=\psi.\\
\end{cases} 
\end{align}
Abusing notation, we still denote the domain $\pi^{-1}(\cU)$ by $\cU$.
Let $\theta$ be a fixed smooth K\"ahler form on $Y$ and we consider the following perturbed family of complex Monge-Amp\`ere equations on $\pi^{-1}(\cU)$ for $s \in (0, 1)$, 

\begin{align}\label{scmp}
\begin{cases}
&(s\theta + \ddbar \varphi_s)^n = e^{\varphi_s} \frac{\prod_i(|\sigma_{E_i}|^2_{h_{E_i}}+ s)^{a_i}}  {\prod_j(|\sigma_{F_j}|^2_{h_{F_j}}+ s)^{b_j}}\Omega_Y, \\
&\varphi_{s_{|\partial\cU}}=\psi.
\end{cases}
\end{align}
 For each $s>0$, we shall first get solution $\varphi_s$ for equation \eqref{scmp}. When $s=0$, equation \eqref{scmp} coincides with equation \eqref{ocmp}.
Next we want to use pluripotential theory to get uniform $C^0$ estimate with barrier of $\varphi_s$. 
Similar $C^0$ estimate of degenerate Monge-Amp\`ere equations have been obtained in different settings  such as  on  unit ball in \cite{Kol1}, on singular Stein domain in \cite{GGZ} and on singular variety with Klt singularity in \cite{EGZ,D-P,Z}. The main differences of our geometric domain with previous setting is the following: we consider the log canonical singularity, which means we do not have  $L^p (p>1)$ integrability of right hand side of equations \eqref{scmp}, hence we do not have uniform boundedness  control of $\varphi_s$, which will also cause extra trouble for high order derivative estimates. 
%Secondly, our geometric domain $\Omega$ is not globally strongly pseudoconvex after we blow up the isolated log canonical singularity. 

Firstly, let us recall that our domain $\cU$ is defined by $\cU:=\{\rho<a\}$ in formula $\ref{domain1}$, where $\rho$ is smooth on $\mathbb C^N$. Hence $\pi^*\rho$ is smooth on $\pi^{-1}\cU$, where $\pi$ is a resolution of singularity. Now choose an arbitrary smooth extension $\psi_1$ of $\psi$, which is supported on a neighborhood of $\partial\cU.$ 
 Hence $\pi^*\rho$ is smooth on $\pi^{-1}\cU$. By choosing A  large enough, we define a semipositive (1,1) form $\omega$ on $\cU$ 
\begin{equation}\label{aux}
 \omega:=A\ddbar (\rho-a)+\ddbar\psi_1.\end{equation}  
 Now
$$s\theta+\ddbar\varphi_s=\omega+s\theta+\ddbar(\varphi_s-A(\rho-a)-\psi_1)$$

Let $\phi_s:=\varphi_s-A(\rho-a)-\psi_1,M=A(\rho-a)+\psi_1$, then we can rewrite equations  \eqref{scmp} as a new family of equations with $\phi_s$ as unknown functions and with zero Dirichlet boundary.
\begin{align}\label{scmp11}
\begin{cases}
&(\omega + s\theta + \ddbar \varphi_s)^n =  \frac{e^{\varphi_s+M}\prod_i(|\sigma_{E_i}|^2_{h_{E_i}}+ s)^{a_i}}  {\prod_j(|\sigma_{F_j}|^2_{h_{F_j}}+ s)^{b_j}}\Omega_Y, \\
&\varphi_{s_{|\partial\cU}}=0.
\end{cases}
\end{align}
The above equations hold for $\phi_s$, we abuse notation and still denote the unknown functions by
$\varphi_s$. Notice all the estimate we get for $\phi_s$ also holds for $\varphi_s$ due to the fact that $M$ is independent of $s$.
From now on, we will focus on equations \eqref{scmp11}.

\begin{lemma} For each fixed $s>0$, there exists a unique smooth solution $\varphi_s$ to equation \eqref{scmp11}.
\end{lemma}
\begin{proof} Recall that $\partial\cU:=\{\rho=a\}$, hence $\phi:=\rho-a$ is equal to $0$ on $\partial\cU$. Now one subsolution of equation \eqref{scmp11} will be $A\phi$ by choosing A sufficient large. it is well-known that existence of subsolution implies existence of solution, see for example \cite[Theorem A]{Bou}.
 \end{proof}

 Next, we state the Kodaira Lemma which is very useful in the estimates that follows.
 \begin{lemma}\cite[Lemma 2.62]{KM}\label{lem:kodaira} Let $\pi:Y\to X$ be log resolution of singularity, where $X$ is $\mathbb Q$ factorial. Also let $\omega$ be a K\"ahler form on $X$ and $\theta$ be a K\"ahler form on $Y$.  
 Then there exists a simple normal crossing divisor $D=\sum D_i$ supporting on the exceptional locus,  hermitian metrics $h_{D_i}$ on the line bundle associated to $D_i$ and a sequence of constants $\alpha_i>0,s_0>0$ such that  $$\pi^*\omega + s\sum_i\alpha_i\ddbar\log h_{D_i} > 0,$$ for all $0<s\leq s_0$. By adjusting the coefficients $\alpha_i$ of $D_i$ we may assume that $s_0=1$ and that   there exists a constant $\be>0$ such that $$\pi^*\omega + \sum_i\ddbar\alpha_i\log h_{D_i} > \be\theta.$$ 
For simplicity, we also define
\begin{equation}\label{Kod}
 |\sigma_D|^2_{h_D}:=\prod_i|D_i|^{2\alpha_i}_{h_{D_i}}, \ddbar\log h_D:=\sum_i\alpha_i\ddbar\log h_{D_i}.
 \end{equation}
\end{lemma}
 In the next two sections we obtain uniform estimates for  $\varphi_s$, completing the proof of Theorem \ref{main1}.
 
 \subsection{The $C^0$-estimate.}

\begin{proposition}\label{prop:c0} Let $|\sigma_D|^2_{h_D}$ be defined as in (\ref{Kod}).
For any $\epsilon>0$,  there exists a constant $C_\epsilon>0$ such that   $$C \geq \varphi_s \geq -(2n+\epsilon)\log \Big(-\log|\sigma_D|_{h_D}^2\Big) - C_\epsilon.$$
 \end{proposition}
 \begin{proof}
 First note that all $\varphi_{ s}$ are $C\theta$- PSH for some fixed large constant $C$ independent of $s$ and $\varphi_{ s}=0$ on $\partial\cU$, this implies that $$\Delta_\theta\varphi_s\geq -C.$$ The upper bound then follows easily by comparing $\varphi_{s}$ to the solution to the following Dirichlet problem:
\begin{align*}
\Delta_\theta u&=-C,\\
u&=0 \text{ on } \partial\cU.
\end{align*}
For the lower bound we will use a barrier function inspired by  \cite[Proposition 3.1]{SSW} (See also \cite{DGG} for an interesting $C^0$ estimate for singular K\"ahler-Einstein metrics in families). From now on, denote $\epsilon>0$  a small but fixed constant and $C_\epsilon$  a constant depending on $\epsilon$ which might change from line to line. Let $D$ be as in  \ref{Kod}. We scale $h_D$ so that $|\sigma_D|_{h_D}^2 < e^{\frac{-3n}{\epsilon}}$. Note that this has no effect on $\ddbar\log h_D$. Now let $H$ be a function defined on $(-\infty,0]$ satisfying $ H'>0,H''>0$ to be determined. For each fixed divisor $F_j$ and some constant $0<\delta<1$, consider the following barrier function $H(u^j_\epsilon)$, where
$$u^j_\epsilon:=\epsilon\log|\sigma_D|^2-\epsilon^2(-\log|F_j|^2_{h_{F_j}})^{1-\delta}.$$
For $\epsilon>0$ sufficiently small, direct calculation shows that
\begin{align*}\label{dd}
    \omega+\ddbar u^j_\epsilon&=\omega-\epsilon Ric(D)-\epsilon^2(1-\delta)(-\log|F_j|^2_{h_{F_j}})^{-\delta}\ddbar(-\log|F_j|^2_{h_{F_j}})\\
    &+\epsilon^2\delta(1-\delta)(-\log|F_j|^2_{h_{F_j}})^{-1-\delta}\sqrt{-1}\partial\log|F_j|^2_{h_{F_j}}\wedge\bar\partial\log|F_j|^2_{h_{F_j}}\\
\end{align*}
Notice that $-\log|F_j|^2_{h_{F_j}}>1$, therefore (possibly change $\theta$) $$\epsilon^2(1-\delta)(-\log|F_j|^2_{h_{F_j}})^{-\delta}\ddbar(-\log|F_j|^2_{h_{F_j}})\leq \epsilon^2\theta.$$ It follows from Kodaira's Lemma that for $\epsilon$ sufficiently small, one has
\begin{equation}\label{1}\omega-\epsilon Ric(D)-\epsilon^2(1-\delta)(-\log|F_j|^2_{h_{F_j}})^{-\delta}\ddbar(-\log|F_j|^2_{h_{F_j}})\geq \frac{\epsilon}{2}\theta.\end{equation}
If we  write $|F_j|^2=|z_j|^2e^{\varphi_j}$ locally near divisor $F_j$, we have
\begin{equation}\label{2}\sqrt{-1}\partial\log|F_j|^2_{h_{F_j}}\wedge\bar\partial\log|F_j|^2_{h_{F_j}}\geq \frac{\sqrt{-1}dz_j\wedge d\bar z_j}{2|z_j|^2}-\theta\end{equation}

Combine inequalities (\ref{1},\ref{2}) and further assume that $H$ can be choosen such that $H'(u^j_\epsilon)<1$, we have 
\begin{align*}
\omega+s\theta+ \ddbar H(u^j_\epsilon) &\geq \omega + H'\ddbar u^j_\epsilon \\
 &\geq H'(\omega + \ddbar u^j_\epsilon)\\
 &\geq H'\Big(\frac{\epsilon}{4}\theta+\epsilon^2\delta(1-\delta)(-\log|F_j|^2_{h_{F_j}})^{-1-\delta}\frac{\sqrt{-1}dz_j\wedge d\bar z_j}{|z_j|^2}\Big)
\end{align*}

Now define $\hat H=\frac{\sum_j H(u^j_\epsilon)}{m}$ (recall that $m$ is the number of components of $F=\sum F_j$) and $\chi_s=\omega+s\theta+\ddbar \hat H$, direct calculation shows that \begin{equation}\label{3}\chi_s^n\geq C_{\delta,\epsilon,n}(H')^n\frac{dz_1\wedge d\bar z_1\cdots dz_n\wedge d\bar z_n}{\Pi_j \Big(|z_j|^2(-\log|z_j|)^{1+\delta}\Big)}\end{equation}
Now we define $H(x) = -B\log(-x)$ for some $3n>B>0$ to be chosen later, and consider $$\phi_s := \varphi_s - \hat H.$$ Then $\phi_s$ solves the following equation\begin{equation} (\chi_{s}+\ddbar\phi_s)^n = e^{\phi_s+\hat H+M}\frac{\prod_i(|\sigma_{E_i}|^2_{h_{E_i}}+ s)^{a_i}}  {\prod_j(|\sigma_{F_j}|^2_{h_{F_j}}+ s)^{b_j}}\Omega_Y.\end{equation} Let $q_{s}$ be the minimum point of $\phi_{s}$. Without loss of generality we may assume that $q_s$ lies in the interior. Since $\phi_s(x)\rightarrow +\infty$ as $x\rightarrow \mathrm{Supp}(D)$, we may also assume that $q_s\notin \mathrm{Supp}(D)$. Then by the maximum principle, at $q_s$ we have the following estimate 

\begin{equation}\label{4}\chi_{s}^n\leq e^{\phi_s+\hat H+M}\frac{\prod_i(|\sigma_{E_i}|^2_{h_{E_i}}+ s)^{a_i}}  {\prod_j(|\sigma_{F_j}|^2_{h_{F_j}}+ s)^{b_j}}\Omega_Y \leq  C\frac{e^{\phi_s+\hat H}}{\prod_j(|\sigma_{F_j}|^2_{h_{F_j}}+ s)^{b_j}}\theta^n,\end{equation}for some constant $C$ independent of $s$. 
Now locally near $q_s$, $\sigma_D = z_1^{\al_1}\cdots z_k^{\al_n}$ where $(z_1,\cdots,z_n)$ are complex coordinates and $\al_i > 0$ for $i=1,\cdots,k$. Combine inequalities (\ref{3},\ref{4}), we have
$$e^{\phi_s}\geq C_{\delta,\epsilon,n}e^{-\hat H}\frac{(H')^n}{{\Pi_j (-\log|z_j|)^{1+\delta}}}\geq C_{\delta,\epsilon,B}\frac{(-\log|\sigma_D|^2)^{B-n}}{{\Pi_j (-\log|z_j|)^{1+\delta}}},$$
where we have used $b_j\leq 1$ and the following two inequalities 
\begin{equation}
e^{-\hat H}\geq B\log(-\log|\sigma_D|^2_{h_D}), (H')^n\geq (-2B\epsilon\log|\sigma_D|^2_{h_D})^{-n}.\end{equation}
Choosing $\delta<<\epsilon$ and $B=2n+\epsilon$, we have $e^{\phi_s}\geq C_{\epsilon}.$ Therefore
\begin{equation}\varphi_s\geq \hat H-C_\epsilon=-\sum_j\frac{(2n+\epsilon)}{m}\log(-u^j_\epsilon)-C_\epsilon\geq -(2n+\epsilon)\log(-\log|\sigma_D|^2)-C_\epsilon.\end{equation}

 \end{proof}

%\begin{remark}
% It is clear from the proof that if $F$ has at most $k$ simple normal crossing components at any point in the support, then the coefficient in front of the  $\log \Big(-\log|\sigma_D|_{h_D}^2\Big)$ term can be improved to $n+k$. %This is optimal \textcolor{red}{why??}.
%\end{remark}

\begin{corollary}\label{cor:finitevol} There exists a constant $C$ independent of $s$ such that $$\int_{\cU\setminus\mathrm{Supp}(D)}(\omega+s\theta + \ddbar \varphi_s)^n\leq C.$$

\end{corollary}
\begin{proof}
Recall that $h_D$ is chosen such that $|\sigma_D|_{h_D}^2 < e^{-4n}$.  Let $f= -3n \log (-\log |\sigma_D|^2_{h_D}) + A$ for some sufficiently large $A>0$ so that 

$$f \geq \varphi_s $$ on $\partial \cU$ for all $s\in (0,1)$. From the calculations above one can see that $f\in \mathrm{PSH}(\omega+s\theta)$.  For $\vep>0$, we let $$\varphi_{s,\vep} = \varphi_s + \epsilon \log |\sigma_{D}|_{h_{D}}^2.$$  By Proposition \ref{prop:c0} we have that $\varphi_{s,\vep}<f$ near $\mathrm{Supp}(D)$ and on $\partial\cU$. In particular, this implies that the set $\{\varphi_{s,\epsilon}>f\}$ will be supported on a relatively compact set contained in $(\cU\setminus Supp(D))$. So by the comparison principle (cf.\cite[Lemma 3.4]{CKZ}) and the fact that $s\theta+\epsilon\ddbar\log h_D \geq 0$ for sufficiently small $\epsilon$, we have, $$\int_{\varphi_{s, \epsilon} > f } (\omega + s\theta + \ddbar \varphi_{s})^n \leq \int_{\varphi_{s, \epsilon} > f } (\omega + 2s\theta + \ddbar \varphi_{s, \epsilon})^n \leq \int_{\varphi_{s, \epsilon} > f } (\omega + 2s\theta + \ddbar f)^n.$$ If we note that as $\vep \searrow 0$, the open sets $W_\vep = \{\varphi_{s,\vep} \geq f\}$ increase to $W = \{\varphi_s>f\}$. Note that $\mathrm{Supp}(D)\subset W$ by Proposition \ref{prop:c0}. We shall see in the next section, that on compact subsets of $\cU\setminus\mathrm{Supp}(D)$ (for instance on $\cU\setminus W$), one has uniform second derivative bounds on $\varphi_s$ (cf.\ Lemma \ref{lem:uniform-smooth-estimates}). So the Corollary immediately follows as long as we can prove that $$ \int_{\cU\setminus \mathrm{Supp}(D)} (\omega + s\theta + \ddbar f)^n<C$$ for some uniform constant $C$. Similar to the calculations above, we have that 
\begin{align*}
\omega+s\theta+\ddbar f & \leq C\theta +  \frac{2n}{(-\log|\sigma_D|_{h_D}^2)^2} \sum_{j=1}^k \al_j\frac{\sqrt{-1}dz_j\wedge d\bar z_j}{|z_j|^2}\\
&\leq C\theta +  C \sum_{j=1}^k \al_j\frac{\sqrt{-1}dz_j\wedge d\bar z_j}{(-\log|z_j|^2)^2|z_j|^2}\
\end{align*}
for some uniform constant $C$. In the second line we also used the fact that $-\log|\sigma_D|_{h_D}^2 \geq -c\log|z_j|^2$ for all $j$ where $c = \min \al_j$. The required estimate then follows from the binomial theorem and the elementary observation that 
\begin{align*}
\int_{|z|<1/2}\frac{\sqrt{-1}dz\wedge d\bar z}{(-\log|z|)^2|z|^2} = \int_0^{1/2}\frac{1}{(-\log r)^2}\cdot\frac{dr}{r} = \int_{\log 2}^{\infty} \frac{ds}{s^2} < \infty. 
\end{align*}
%Note that the observation above is simply the well known fact that the Poincar\'e metric on the punctured disc has finite volume. 

 \end{proof}

\subsection{Further estimates and proof of Theorem \ref{main1}.}~\\

We first  prove the boundary $C^1$ estimate. We denote the covariant derivative of $\theta$ by $\nabla$. We also let $D$ be the effective divisor from Lemma \ref{lem:kodaira}, such that $\mathrm{Supp}(D) = \mathrm{Supp}(F)\cup \mathrm{Supp}(E)$ and 
\begin{equation}\label{eq:kodaira}
\omega+\ddbar\log h_D > \be\theta
\end{equation} for some $\be>0$. Note that by the support condition, if $\sigma_D$ is a defining section of $D$, then there exists a uniform constant $C$ such that 
\begin{equation}\label{eq:estforsections}
|\sigma_F|^{-2},|\sigma_E|_{h_E}^{-2}\leq C|\sigma_D|^{-2l},
\end{equation}  for some $l\in \mathbb{N}$.
\begin{lemma}\label{BC1} There exists a constant $C>0$ such that for all $0<s<1$,   $$|\na\varphi_{ s}|_{\partial\cU}\le C.$$ 
\end{lemma}
\begin{proof}
 As in the proof of Proposition \ref{prop:c0}, let $u$ solve the Dirichlet problem \begin{align*}
\Delta u=-A,\\
u=0 \text{ on } \partial\cU.
\end{align*}
Then if we choose $A$  to be large enough, we have that $\varphi_s\leq u$.  Next, again by Proposition \ref{prop:c0}, $$\varphi_{ s}\text{ and }\frac{e^{\varphi_s}\prod_i(|\sigma_{E_i}|^2_{h_{E_i}}+ s)^{a_i}}  {\prod_j(|\sigma_{F_j}|^2_{h_{F_j}}+ s)^{b_j}} \frac{\Omega_Y}{(\omega+s\theta)^n}$$ are uniformly bounded in the neighbourhood of the boundary
 $\partial\cU$. Let $U_\vep := \{\rho > a-\vep\}$. We fix a small $\vep$ such that $U_\vep$ has smooth boundary consisting of two components $\partial \cU$ and $\{\rho = a-\vep\}$. Then we can choose $b>>1$ such that $$(\omega+ s\theta+ \ddbar [b(\rho-a)])^n >  \frac{e^{\varphi_s+M}\prod_i(|\sigma_{E_i}|^2_{h_{E_i}}+ s)^{a_i}}  {\prod_j(|\sigma_{F_j}|^2_{h_{F_j}}+ s)^{b_j}}\Omega_Y = (\omega+s\theta + \ddbar\varphi_s)^n.$$ 
 
  %Now we show that the maximum of $b(\rho-a)-\varphi_s$ must be achieved at the boundary $\partial U_\epsilon$. If not, then at some interior point $p$ of $U_\epsilon$, we have the following contradiction
  %\begin{align*}0&<\det((f+b(\rho-a))_{i\bar j})-\det((f+\varphi_s)_{i\bar j})\\&=\sum_{ij} \Big(\int_0^1 B^{i\bar j}(t)dt\Big)\Big(b\rho-ba-\varphi_s\Big)_{i\bar j}\leq 0
  %\end{align*}
 %where $f$ is the K\"aher potential of $\omega+s\theta$ near $p$ and $B^{i\bar j}(t)$ is the cofactor of $$\Bigg(t(f+b\rho-ba)+(1-t)(f+\varphi_s)\Bigg)_{i\bar j}.$$
 On the other hand we also have 
 \begin{align*}
 b(\rho-a)\Big|_{\partial \cU}  = 0 = \varphi_s\Big|_{\partial \cU}\\
 b(\rho-a)\Big|_{\rho = a-\vep}  = -b\vep < \varphi_s\Big|_{\rho = a-\vep},
 \end{align*}
 if we pick $b>>1$. Then by the maximum principle, and the upper bound above, $$b(\rho - a)\leq \varphi_s \leq u.$$ But then it is easy to see that for any $x\in \partial \cU$, $$|\nabla \varphi_s|(x) \leq \max\Big(|\nabla b(\rho-a)|(x),|\nabla u|(x) \Big),$$ and the boundary $C^1$-estimate follows.%Let $U\subset\cU$ be a neighbourhood of the boundary $\partial \cU$, say $U = \{\rho\geq a-\vep\}$ for some $\vep<<1$. We choose $\vep$ so that $\{\rho = a-\vep\}$ is a smooth manifold.
 \end{proof}

Next, we prove the global $C^1$ estimate with suitable barrier function.  Such gradient estimate without barrier is  firstly studied in \cite{Blo} for standard non-degenerate Monge-Amp\`ere equation where $\varphi$ is bounded and later improved in \cite{PS} when the potential only has an upper bound in a different geometric setting. 
\begin{proposition}\label{GC1}
Let $\varphi_{s}$ be the solution of equation  \eqref{scmp11}, There exist $N,C>0$ such that for all $0<s<1$,   $$|\na\varphi_{ s}|^2|\sigma_D|^{2N}_{h_D}\le C.$$  
\end{proposition}
\begin{proof} 
 Once again we fix a constant $\be>0$ such that $\omega+\ddbar\log {h_D}>\be\theta$ and rewrite the equation \eqref{scmp11} as
 $$\left(\omega + \ddbar\log h_D + s\theta+ \ddbar\phi_{ s} \right)^n = \frac{e^{\varphi_s+M}\prod_i(|\sigma_{E_i}|^2_{h_{E_i}}+ s)^{a_i}}  {\prod_j(|\sigma_{F_j}|^2_{h_{F_j}}+ s)^{b_j}} \Omega_Y.$$
where $$\phi_{ s}=\varphi_{ s}-\log |\sigma_D|^2_{h_D}.$$
It follows from Proposition \ref{prop:c0} that $\phi_s>-C$ for some uniform constant $C>0$. Note that our reference metrics $ \omega +\ddbar\log h_D + s\theta$ in the above equation are uniformly non degenerate as $s\rightarrow 0$. By abusing notation, in the rest of proof, we will use $\theta$ to denote the K\"ahler form $\omega+\ddbar\log h_D+s\theta$. Define
$$F:=\frac{e^{\varphi_s+M}|\sigma_D|^2_{h_D}\prod_i(|\sigma_{E_i}|^2_{h_{E_i}}+ s)^{a_i}}  {\prod_j(|\sigma_{F_j}|^2_{h_{F_j}}+ s)^{b_j}}\frac { \Omega_Y}{( \omega + \ddbar\log h_D + s\theta)^n}$$
 
For the rest of the proof, for convenience we drop the sub-script $s$ from the notation ie. we denote $\phi_s$ by simply $\phi$. We define
$$H=\log|\nabla\phi|_\theta^2+\log |\sigma_ D|_{h_D}^{2k}-\gamma(\phi),$$ 
where $k$ is a constant, $\gamma$ is a one variable monotone increasing  function to be determined. 
 Since the leading term of our function $\gamma(x)$ will be chosen as $Bx$, and $\phi$ blows up in the rate of $-\log |\sigma_D|^2_{h_D}$, in particular $-\gamma(\phi)$ has upper bound. Now $\log|\nabla\phi|_\theta< C-2\log|\sigma_D|_{h_D}$,  this implies that 
 when 
$k>2$, $H$ has a maximum in $\cU\setminus \mathrm{Supp} (D)$.
Direct computation (cf. \cite[Page 21]{PSS1}) shows that
\begin{equation}\label{PS}
\Delta'\log \,|\nabla\phi|_\theta^2
\ge
\frac{2 {\Re}\langle\nabla\log F, \nabla \phi\rangle_\theta
}{ |\nabla\phi|_\theta^2}
-\Lambda \,tr_{\omega'}\theta
+
\frac{|\nabla\nabla\phi|_{\theta,\omega'}^2+|\bar\nabla\nabla\phi|_{\theta,\omega'}^2}{
|\nabla\phi|_\theta^2}-
\frac{|\nabla|\nabla\phi|_\theta^2|_{\omega'}^2}{ |\nabla\phi|_\theta^4}
\end{equation}
where $\Delta'$ is taken with respect to metric $\omega'=\omega +\ddbar\log h_D + s\theta+ \ddbar\phi_s$, $\Lambda$ is the bound of bisectional curvature of metric $\theta$ and $|\cdot |_{\theta,\omega'}$ means that the norm of the two tensor is taken with respect to $\theta$ on the first entry and $\omega'$ on the second entry. When $|\nabla\phi|_\theta>1$, by Cauchy-Schwarz inequality,  
\begin{equation}\label{3.9}
|\frac{2 {\Re}\langle\nabla\log F, \nabla\phi\rangle_{\theta}
}{ |\nabla\phi|_\theta^2}|\le C+C\sum_i|\sigma_{E_i}|^{-2}_{h_{E_i}}+\sum_j|\sigma_{F_j}|^{-2}_{h_{F_j}}+|\sigma_D|^{-2}_{h_D}\le C|\sigma_D|^{-2l}_{h_D},
\end{equation}
for some $l>0$. As observed before, the final equality follows from the fact  that $\mathrm{Supp}(D) = \mathrm{Supp}(E)\cup \mathrm{Supp}(F)$. Next, by \cite[Lemma 13]{PSS1}, we have

\begin{equation}\label{3.10}\frac{|\na\na \phi|_{\theta,\omega'}^2+|\na\bar\na \phi|^2_{\theta,\omega'}}{ |\na\phi|_\theta^2}
- \frac{|\na |\na\phi|_\theta^2|_{\omega'}^2}{ |\na\phi|_\theta^4}  
\ge 2\Re\langle \frac{\na|\na\phi|_\theta^2}{ |\na\phi|_\theta^2}, \frac{\na \phi}{ |\na\phi|_\theta^2}\rangle_{\omega'} - 
2\Re\langle \frac{\na|\na\phi|_\theta^2}{ |\na\phi|_\theta^2}, \frac{\na \phi}{ |\na\phi|_\theta^2}\rangle_\theta 
\end{equation}
At the maximum of $H$, we have 
\begin{align*}\na\log|\na\phi|^2+\na\log|\sigma_ D|_{h_D}^{2k}-\gamma'\na\phi=0\end{align*}
Hence we have 
\begin{equation}\label{3.11}
\begin{split}&2\Re\langle \frac{\na|\na\phi|_\theta^2}{ |\na\phi|_\theta^2}, \frac{\na \phi}{ |\na\phi|_\theta^2}\rangle_{\omega'} - 
2\Re\langle \frac{\na|\na\phi|_\theta^2}{ |\na\phi|_\theta^2}, \frac{\na \phi}{ |\na\phi|_\theta^2}\rangle_\theta\\
&=2k\Re\langle -\na\log|\sigma _D|_{h_D}^2+\gamma'\na\phi,\frac{\na\phi}{|\na\phi|^2_\theta} \rangle_{\omega'}-
2k\Re\langle -\na\log|\sigma_ D|_{h_D}^2+\gamma'\na\phi,\frac{\na\phi}{|\na\phi|^2_\theta} \rangle_\theta\\
&\ge 2k\Re\langle -\frac{\na|\sigma _D|_{h_D}^2}{{|\sigma_D|}^2_{h_D}},\frac{\na\phi}{|\na\phi|^2_\theta} \rangle_{\omega'}+2k\gamma'\frac{|\nabla\phi|^2_{\omega'}}{|\nabla\phi|^2_\theta}+
2k\Re\langle \frac{\na|\sigma_ D|_{h_D}^2}{|\sigma_D|^2_{h_D}},\frac{\na\phi}{|\na\phi|^2_\theta} \rangle_\theta-2k\gamma'
\end{split}
\end{equation}
We will choose $\gamma$ as a monotone increasing function, so we can drop the term $2\gamma'\frac{|\nabla\phi|^2_{\omega'}}{|\nabla\phi|^2_\theta}$. At the maximum of $H$, we can assume $|\sigma_D|^{2k}_{h_D}|\na\phi|^2_\theta\ge 1$ otherwise we are done. Choosing $k\ge 2$, we have the following two estimates
\begin{equation}\label{3.12}
\begin{split}|2\Re\langle -\frac{\na|\sigma _D|_{h_D}^2}{|\sigma_D|_{h_D}^2},\frac{\na\phi}{|\na\phi|^2_\theta} \rangle_{\omega'}|&\le 2|\Re\langle \na|\sigma _D|_{h_D}^2,\frac{\na\phi}{|\na\phi|_\theta}\rangle_{\omega'}|\\
&\le
|\na|\sigma_D|^2_{h_D}|^2_{\omega'}+\frac{|\sigma_D|^2_{h_D}|\na\phi|^2_{\omega'}}{|\sigma_D|^2_{h_D}|\na\phi|^2_{\theta}}\\
&\le 
C|\na|\sigma_D|^2_{h_D}|^2_{\theta}\,tr _{\omega'}\theta+|\sigma_D|^2_{h_D}|\na\phi|^2_{\omega'}.\\
|2\Re\langle -\frac{\na|\sigma _D|_{h_D}^2}{{|\sigma_D|}^2_{h_D}},\frac{\na\phi}{|\na\phi|^2_\theta} \rangle_{\theta}|&\le 2|\Re\langle \na|\sigma _D|^2_{h_D},\frac{\na\phi}{|\na\phi|_\theta}\rangle_{\theta}|\le C.
\end{split}\end{equation}
On the other hand, one has
\begin{equation}\label{2.19}-\Delta'\gamma(\phi)=-\gamma'\Delta'\phi-\gamma''|\na\phi|^2_{\omega'}
=\gamma'tr_{\omega'}\theta-n\gamma'-\gamma''|\na\phi|^2_{\omega'},\Delta'\log|\sigma_D|^k_{h_D}\le C\,tr_{\omega'}\theta.\end{equation}
Combine estimate \eqref{2.19} with preceding estimates \eqref{3.9}, \eqref{3.10}, \eqref {3.11}, \ref {3.12}), we have 
\begin{equation}\label{3.13}
\Delta' H\ge (\gamma'-\Lambda-C)tr_{\omega'}\theta-(n+2)\gamma'-
(\gamma''+|\sigma_D|^2_{h_D})|\na\phi|^2_{\omega'}
-C|\sigma_D|^{-2l}_{h_D}.
 \end{equation}
Recall that $\phi>-C'$, now we construct our function $\gamma$ as
$$\gamma(x)=(\Lambda+C+1)x-\frac{E}{x+C'+1}$$
where $E$ is a constant to be determined.
Then by \eqref{3.13} we have
\begin{align*}
\Delta' H\ge tr_{\omega'}\theta-(n+2)(C+1+\Lambda)-C|\sigma_D|^{-2l}_{h_D}+(\frac{2E}{(\phi+C'+1)^3}-|\sigma_D|^{2}_{h_D})|\na\phi|^2_{\omega'}
\end{align*}
Noticing that $\phi\le C-2\log|\sigma_D|^2_{h_D}$,  we can  assume that
$$(\frac{2E}{(\phi+C'+1)^3}-|\sigma_D|^2_{h_D})\ge |\sigma_D|^2_{h_D}$$
by choosing suitable large constant $E$ depending on $C,C'$. Now we conclude that at the maximum of $H$,  we have

\begin{align*}tr_{\omega'}\theta\le (|\sigma_D|_{h_D})^{-2l}, ~|\na\phi|^2_{\omega'}\le C(|\sigma_D|_{h_D})^{-2(l+1)}.\end{align*}
Hence$$|\na\phi|^2_\theta\le C(|\sigma_D|_{h_D})^{-2(2l+1)}.$$ 
Choosing $k=2l+1$, we have  $H_{max}\le C$, and letting $N:=k+(C+1+\Lambda)$, it follows from the definition of $\gamma(\phi)$ that
\begin{align*}|\na\phi|_\theta^2|\sigma_D|_{h_D}^{2N}\le C.\end{align*}
\end{proof}

\begin{lemma}\label{BC2}
Let $\varphi_{ s}$ be the solution of equation  \eqref{scmp11}, then there exist constant $C>0$ such that for all $0<s<1$,    $$|\na_\theta^2\varphi_{ s}|_{\partial\cU}\le C.$$
\end{lemma}
\begin{proof}
 Notice that our boundary is strictly pseudoconvex, and all data in the equation \eqref{scmp11} is uniformly bounded near the boundary, so the second order estimates on the  boundary  follow directly from the arguments in \cite[Section 1.3]{CKNS}. 
\end{proof}

Next, we will prove  second order estimates with bounds from suitable barrier functions. %There exists an effective Cartier divisor $D$ on $Y$ such that for any sufficiently small $s>0$,

\begin{lemma}\label{2ndem} There exist constants $N,C>0$ such that for all $0<s<1$,  
\begin{equation}
\sup_{\cU} \left(  |\sigma_D|^{N}_{h_D} \right)  |\Delta_\theta \varphi_{ s} | \leq C,
\end{equation}
where $\Delta_\theta$ is the Laplace operator with respect to the K\"ahler metric $\theta$.
\end{lemma}
\begin{proof} 
We remark that the constant $C$ in the proof might change from line to line and it depends on $\theta,\mathcal U$ but does not depend on $s$.
Let $\omega' =  \omega + s\theta + \ddbar \varphi_{s}$. Then we consider the quantity $$H= \log tr_\theta (\omega') - B\varphi_{ s} +B \log |\sigma_D|^2_{h_D}$$
for some large constant $B>2$. By the $C^0$ estimate of $\varphi_{ s}$, $H$ is bounded above in $\cU$. Standard calculations (cf. \cite[Lemma 3.7]{Gabor}) show that 
\begin{equation}\label{eq:aubin-yau-c2}
\Delta ' \log tr_{\theta}{\omega'} \geq -Ctr_{\omega'}{\theta} - \frac{tr_{\theta}Ric(\omega')}{tr_{\theta}(\omega')}. 
\end{equation} 
where $C$ depends on bisectional curvature of $\theta$.
From equations \eqref{scmp11}, \eqref{eq:estforsections} and the elementary observation that $$\ddbar\log(f+s) = \frac{f}{f+s}\ddbar\log f + s\frac{\sqrt{-1}\partial f\wedge \overline{\partial} f}{f(f+s)}$$ holds for any smooth non-negative function $f$, it is easy to see that $$-tr_{\theta}Ric(\omega') \geq \frac{-C}{|\sigma_D|_{h_D}^{2l}},$$ for some constant $C>0$ independent of $s$.  Together with \eqref{eq:aubin-yau-c2} and our choice of $D$ (cf.\ \eqref{eq:kodaira})  we see that 
\begin{align*}
\Delta'H &\geq -Ctr_{\omega'}{\theta}  - \frac{C}{|\sigma_D|_{h_D}^{2l}tr_{\theta}(\omega')} + Btr_{\omega'}(\ddbar \log h_D - \ddbar\varphi_s)\\
&\geq (B\be - C)tr_{\omega'}{\theta} - \frac{C}{|\sigma_D|_{h_D}^{2l}tr_{\theta}(\omega')} - Bn\\
&\geq tr_{\omega'}{\theta} - \frac{C}{|\sigma_D|_{h_D}^{2l}tr_{\theta}(\omega')} - Bn \hspace{0.5in} (\text{if } B>>1)\\ 
&\geq (tr_\theta (\omega'))^{\frac{1}{n-1}}(\frac{\theta^n}{\omega'^n})^{\frac{1}{n-1}} - \frac{C}{|\sigma_D|_{h_D}^{2l}tr_{\theta}(\omega')} - Bn\\
&\geq  (tr_\theta (\omega'))^{\frac{1}{n-1}}|\sigma_D|_{h_D}^{2\al} - \frac{C}{|\sigma_D|_{h_D}^{2l}tr_{\theta}(\omega')} - Bn,
\end{align*}
for some constant $\al$ independent of $s$ and $B>>1$ so that $B\be > C+1$ in line three. 

% if $B$ is large enough, then there exists $C>0$ such that for all $0<s<1$, 
%
%\begin{eqnarray*}
%
%\Delta_{\omega'} H &\geq & \Delta_{\omega'}\log tr_{\theta}(\omega')+2B^2 tr_{\omega'}(\omega' + \ddbar\log h_D)-nB^3  \\
%&\geq & Btr_{\omega'}\theta-C\frac{tr_{\theta}Ric(\omega')}{tr_{\theta}(\omega')}-B^3\\
%
%& \geq&  B(tr_\theta (\omega'))^{\frac{1}{n-1}}(\frac{\theta^n}{\omega'^n})^{\frac{1}{n-1}} - C\frac{1}{tr_\theta\omega' |\sigma_D|^{2\alpha}_{h_D}} -B^3\\
%& \geq&  B(tr_\theta (\omega'))^{\frac{1}{n-1}}|\sigma_D|_{D}^{2\beta} - C\frac{1}{tr_\theta\omega' |\sigma_D|^{2\alpha}_{h_D}} -B^3\\
%
%\end{eqnarray*}
%
%where $\alpha,\beta$ are fixed constants only depending on the coefficient of exceptional divisor in the log resolution of singularity. For the third and fourth inequality,  we use the equation. Noticing by the $C^0$ estimaste of  $\varphi_{ s}$ that for some fixed $\epsilon_0$ to be determined,  $$H<\log(tr_{\theta}(\omega')|\sigma_{D}|^{B^2-B^3\epsilon_0}_{h_D})+C_{\epsilon_0}$$
%
By Lemma \ref{BC2}, it suffices to assume that $ H$ obtains maximum at a point $p\in \cU$. Moreover, since $H$ goes to $-\infty$ on $\mathrm{Supp}(D)$, clearly $p\notin \mathrm{Supp}(D)$.  From the maximum principle it follows that at point $p$,
\begin{equation} \frac{C}{|\sigma_D|_{h_D}^{2l}tr_{\theta}(\omega')}+Bn\geq (tr_\theta (\omega'))^{\frac{1}{n-1}}|\sigma_D|_{h_D}^{2\al}.\end{equation}
We first assume that $|\sigma_D|^{2l}_{h_D}tr_\theta(\omega')\geq 1$ at $p$, then \begin{equation}(tr_\theta (\omega'))^{\frac{1}{n-1}}|\sigma_D|_{h_D}^{2\al}\leq C+Bn.\end{equation}
Notice that the other case is $|\sigma_D|^{2l}_{h_D}tr_\theta(\omega')\leq 1$ at $p$,
if follows that in both cases there is an integer $k$ (depending on $\al,l$ and $n$) such that
\begin{equation}|\sigma_D|_{h_D}^{2k}(tr_\theta\omega')(p)\leq C+Bn.\end{equation} Notice that $\varphi_s\geq \epsilon\log |\sigma_D|^2_{h_D}-C_\epsilon$ for any $\epsilon>0$, so it follows  by choosing a $B>>k$ that 
$H(p)\leq C+Bn$. Now fixing this $B$, we have $H(x)\leq C$ 
for any $x\in \cU\setminus\mathrm{Supp}(D)$. By the definition of $H$, we have
\begin{align*}
|\sigma_D|^{2B}_{h_D}(tr_{\theta}\omega')(x) \leq C.
\end{align*}  Choosing $N=2B$ will finish proof.    %From the equation,  we also have that $|\sigma_D|^N_{h_D}(tr_{\omega'}\theta)(x) \leq C,$ possibly by choosing an even larger $N$. This completes the proof of the Lemma. 
%Also, since our goal is to bound $H$, we may assume that at the maximal point of H, $$(tr_{\theta}(\omega')|\sigma_D|_{h_D}^{2\alpha})>1.$$ Applying the maximum principle, at the maximal point $x_{max}$ of $H$, (by choosing $A> 2\beta(n-1)$)
%
%$$H(x_{max})<tr_\theta(\omega')|\sigma_D|_{h_D}^{A^2}+C  \leq n A^2. $$ 
%
%Hence $\max H\le \max\{C_{\epsilon_0},nA^2\}$. On the other hand, since $\varphi_{q,s,\delta,\epsilon}\leq C$,
%we have $$tr_\theta{\omega'}\le \frac{1}{|\sigma_{D}|_{h_D}^{2A^2}}$$
%This proves the lemma. 

\end{proof}

The following lemma on local higher order regularity of $\varphi_{ s}$ is established by the standard linear elliptic theory after applying Lemma \ref{2ndem} and linearizing the complex Monge-Amp\`ere equation \eqref{scmp11}.

\begin{lemma} \label{lem:uniform-smooth-estimates}
For  any compact  $K \subset\subset (\cU \setminus p) $, and any natural number $k\in \mathbb{N}$, there exists a constant $C=C( k, K)>0$ such that for any $0<s<1$ 
$$ || \varphi_{ s} ||_{C^k(K)} \leq C. $$

\end{lemma}

In our context, we need a lemma for compactness of quasi PSH function in $L^1$ topology. This kind of lemma is standard and  we include it for convenience.
\begin{lemma}\label{com}\cite{GZ} Let $\cU, \theta$ and $\varphi_s$ be as above. Denote $A(\rho-a)+\psi_1+\varphi_s$ by $\psi_s$ (cf. \eqref{aux} for the definition of $A(\rho-a)+\psi_1$). Then $\psi_s$ converge to a PSH function on $\pi^{-1}(\cU)$ by taking a subsequence. 
\end{lemma} 
\begin{proof}Fix a measure $\theta^n$ and a finite open covering $V_i$ for $\pi^{-1}(\mathcal U)$ such that $\theta$ is $\ddbar$ exact on each chart. On each $V_i$,  $\psi_s$ is a $s\theta$- PSH function  satisfying $\|\psi_s\|_{L^1}<C$ and $\psi_s<C$, for some $C$ independent of $s$ (cf. Proposition \ref{prop:c0}). By passing to a subsequence, one has $$\psi_s\xrightarrow{L^1(V_i)} \varphi_i,$$ where  $\varphi_i$ is a $s\theta$- PSH function for any $s>0$ (cf. \cite{Demailly}). Passing $s\to 0$, $\varphi_i$ is PSH on $V_i$. For different charts $V_i,V_j$ with $V_i\cap V_j\neq\emptyset$, we aim to show that $\varphi_i=\varphi_j$. Then we can patch all $\varphi_i$ to get a PSH function on $\pi^{-1}(\mathcal U)$. By taking a further subsequence,  we may assume that $\varphi_i=\varphi_j$ on a full measure subset $W$ of $V_i\cap V_j$. Hence it suffices to show that for any point $q\in V_i\cap V_j$,
$$\varphi_i(q)=\limsup_{x\to q, x\in W}\varphi_i(x).$$
It follows from the upper semicontinuity of $\varphi_i$ that
$$\varphi_i(q)\geq\limsup_{x\to q, x\in W}\varphi_i(x).$$
To get the reversed inequality, we argue by contradiction.
If not, then 
$$\varphi_i(q)>\limsup_{x\to q, x\in W}\varphi_i(x).$$
By the mean value inequality for PSH function, one has that, for some sufficiently small $r>0$, 
$$c_r\int_{B(q,r)\cap W}\varphi_i dV<\varphi_i(q)\leq c_r\int_{B(q,r)}\varphi_i dV,$$
where $B(q,r)$ is a radius $r$ ball in $V_i\cap V_j$ and $c_r$ is a constant depending on $r$.
This is a contradiction by noticing that $\varphi_i<C$ and $B(q,r)\setminus W$ has measure zero.
\end{proof}
Now we proceed to prove our first main Theorem.
\subsection*{Proof of Theorem \ref{main1}}
%\begin{proof} [Proof of Theorem \ref{main1}.] 
By Lemma \ref{com}, for any sequence $s_j\to 0$,  one has 
\begin{equation}
\varphi_{s_j}\xrightarrow{L^1(\pi^{-1}(\cU))} \varphi,
\end{equation}
by taking a subsequence and moreover if we define \begin{equation}\label{def}\varphi_{\KE}:=A(\rho-a)+\psi_1+\varphi,\end{equation} then $\varphi_{\KE}$ is a PSH function on $\pi^{-1}(\cU)$.
It follows from the uniform estimates for $\varphi_{ s}$ away from $\mathrm{Supp}(D)$ that by taking a  further subsequence, one has 
$$\varphi_{s_j}\xrightarrow{C^\infty_{loc}(\pi^{-1}(\cU)\setminus\mathrm{Supp}(D))} \varphi.$$
Clearly, $\varphi$  solves the desired M\"onge-Ampere equation on $(\bar{\cU}\setminus \{p\})$ by design. Moreover, by Proposition \ref{prop:c0} there exists $C>0$  such that 
\begin{equation}\label{C0}
-(2n+\epsilon) \log(-\log |\sigma_D|^2_{h_D}) - C_\epsilon \leq \varphi \leq  C. \end{equation}
With the above sublog $C^0$ estimate of $\varphi$ and the uniqueness Lemma \ref{uniq} below, we have 
$$\varphi_s\xrightarrow{C^\infty_{loc}(\pi^{-1}(\cU)\setminus\mathrm{Supp}(D))} \varphi.$$
when $s\to 0$.
Now $\varphi_\KE$ satisfies properties $(1)-(2)$ in the statement of Theorem \ref{main1} and  $\omega_{\KE}:=\ddbar\varphi_{\KE}$ is  a K\"ahler-Einstein metric on $\cU\setminus\{p\}\cong\pi^{-1}(\cU)\setminus\mathrm{Supp}(D)$. From Corollary \ref{cor:finitevol} and local smooth convergence on $\pi^{-1}(\cU)\setminus\mathrm{Supp}(D)$, we also get by Fatou Lemma that $$\int_{\cU\setminus\{p\}}\omega_\KE^n \leq \lim_{s\rightarrow 0}\int_{\cU\setminus \{p\}}(\omega + s\theta + \ddbar\varphi_s)^n < \infty,$$ and this proves $(3)$.

Part (4) We remark that in \eqref{def}, our function $\varphi_{\KE}$ is PSH function on $\pi^{-1}(\cU)$.
Notice that $\pi^{-1}(\cU)$ is smooth, by the sublog pole estimate of $\varphi_{\KE}$ in Proposition \ref{prop:c0}, it follows  that $\varphi_{\KE}$ has vanishing Lelong number at any point of $\pi^{-1}(\cU)$. 

At last, we show that $\lim_{x\to\mathrm{Supp}( D)}\varphi_{\KE}(x)=-\infty$ for $x\in\pi^{-1}(\cU)$.  We argue by contradiction. If not,  noticing that $\varphi_{\KE}$ is bounded from above, then there is a smooth irreducible component $E$ of $D$ and a sequence of point $y_i\to y\in E$  such that \begin{equation}\label{finite}C\leq\lim_{i\to +\infty}\varphi_{\KE}(y_i)\leq \varphi_{\KE}(y)\end{equation}
for some finite constant $C$. The second inequality is due to the upper semicontinuity of $\varphi_{\KE}$. Now by the definition of PSH function (cf. \cite[Definition 1.4]{Demailly}), $\varphi_{\KE}$ is still a PSH function when restricted to a  smooth component $E$ of $D$ . Indeed, a function with value in $[-\infty,+\infty)$ is PSH if and only if it is upper semicontinuous and it satisfies the mean value inequality when restricted to any complex line. One can check easily that the mean value inequality on complex line and the upper semicontinuity  are preserved  when restricted to complex submanifold. 
%Also for point $q\in E$, it follows from the upper semicontinuity of $\varphi_{KE}$ that 
%\begin{align*}
    %\varphi_{KE}(q)\geq\limsup_{z\in E, z\to q}\varphi_{KE}(z).
%\end{align*}
%The reversed inequality is trivially true. So the upper semicontinuity of $\varphi_{KE}$ is also proved.

Now we have showed that $\varphi_{\KE}|_E$ is a PSH function and moreover it is finite at a point $y\in E$ (cf \eqref{finite}), we claim that $\varphi_{\KE}$ is constant on $E$. Indeed, since $\varphi_{\KE}|_E$ is upper semicontinuous, we may assume $\varphi_{\KE}|_E$ achieve its maximum at some point $x_{max}$.  Then by the mean value inequality for PSH function, it is clear that $\varphi_{\KE}$ is locally constant. Now $E$ is a smooth manifold without boundary, hence  $\varphi_{\KE}|_E$ is constant on $E$. The claim is proved. Next, using the connectedness of $\mathrm{Supp}(D)$ (using $p$ is a normal singularity), $\varphi_{\KE}$ must be constant on $\mathrm{Supp}(D)$. Then we need the following general fact from \cite[Lemma 2.7]{B}: Suppose $\phi$ is a plurisubharmonic function on the unit ball $B \subset \mathbb{C}^n$ such that $$\int_{B} |z_1|^{-2}e^\phi (\sqrt{-1})^ndz_1\wedge d\overline{z_1}\wedge... \wedge dz_n \wedge d\overline{z_n} <\infty, $$ then $\phi$ tends to $-\infty$ near $B \cap \{z_1=0\}$. Now pick a smooth exceptional divisor $F$ with discrepancy $-1$ from $D$ and also pick a generic point $q\in F$ ($q$ is not included in other components of $D$),  then applying the fact  just recalled, one has, $\varphi_{\KE}(x_i)$ tends to $-\infty$ when $x_i\to q$. This is a contradiction since $\varphi_{\KE}$ is a finite constant on $\mathrm{Supp}(D)$.

%Therefore the  function $\varphi$ can uniquely $descend to $(\cU\setminus $p)$, and we obtain a $K\"ahler-Einstein metric %on

Finally, the uniqueness is a consequence of the slightly stronger uniqueness theorem below.

%\end{proof}

\begin{lemma} \label{uniq} There exists a unique smooth PSH function $\varphi \in  C^\infty(\cU\setminus p)$ satisfying

\begin{enumerate}
\item $(\ddbar \varphi)^n = e^\varphi \Omega_X$ on $(\cU\setminus p)$, $\varphi|_{\partial\mathcal U}=\psi$.

\item For any $\epsilon>0$, there exists a constant $C_\epsilon>0$ such that
\begin{align*}\epsilon \log |\sigma_D|^2_{h_D} - C_{\epsilon}  \leq \varphi \leq C , \end{align*}
where $D$ is a SNC divisor supported on the  exceptional locus of a log resolution $\pi$ (cf. Lemma \ref{lem:kodaira} for the definition).

\end{enumerate}

%In particular, $\varphi \in PSH(\cU)$ satisfies all the conditions in Theorem \ref{main1}. 

\end{lemma}

\begin{proof}

Let $\varphi$ be the limit of $\varphi_{ s}$  as $s\rightarrow0$ as above.
Suppose there exists another $\varphi' $ satisfying the conditions $(2)$ in the lemma solving equation \ref{scmp11}, and for any $\epsilon>0$, there exist $C_1>0$ and $C_2= C_2(\epsilon)>0$ such that  
\begin{align*}\epsilon \log |\sigma_D|^2_{h_D} - C_2 \leq \varphi' \leq C_1\end{align*}
We consider the quantity 
\begin{align*}\phi = \varphi-  \varphi' + \delta^3\log |\sigma_D|^2_{h_D}+\delta M,\end{align*}
where recall that $\omega=\ddbar M$ in equation \eqref{scmp11}. 
Then $\phi$ satisfies the following equation on the log resolution $Y$, 

$$\frac{(\omega + \ddbar \varphi'   - \delta\omega + \delta^3 Ric(h_D)+   \ddbar \phi  )^n }{ (\omega + \ddbar \varphi')^n} = e^{\phi-\delta M}    |\sigma_D|^{-2\delta^3}.$$

Noticing that $\varphi, M$ are bounded from above on $(\cU\setminus p)$ and $(\delta^3\log|\sigma_D|^2_{h_D}-\varphi')(x)\to -\infty$ when $x\to \mathrm{supp}(D)$, so $\phi$ obtains its maximum in $(\mathcal U\setminus p)$. It follows from maximum principle that there exists $C>0$ such that
$$\sup_{(\mathcal U\setminus p)} \phi \leq \delta C. $$
Therefore for any point $x\in (\mathcal U\setminus p)$, one has 
$$\varphi-\varphi'\leq \delta C-\delta M-\delta^3\log|\sigma_D|^2_{h_D}.$$ 
Let $\delta \rightarrow 0$,   on ($\mathcal U\setminus p$), one has
$$\varphi \leq \varphi'.$$
Similarly, one has $\varphi \geq \varphi'$ on $(\mathcal U\setminus p)$. Therefore $\varphi=\varphi'$.

%by applying the maximum principle to $$\phi' = \varphi_{ s, \delta, \epsilon} - \varphi' - \delta \log |\sigma_D|^2_{h_D}. $$

\end{proof}
\section{Volume Stability} 
The goal of this section is to prove Theorem \ref{volume} and some consequences. Suppose we have two complete K\"ahler-Einstein metric $\omega_{KE}$ and $\omega'_{KE}$ on $(\cU\setminus p)$,
by K\"ahler-Einstein condition, we can write $\omega_{KE}'=\omega_{KE}+\ddbar\varphi$, where  $$\varphi:=
\log\frac{\omega_{KE}'^n}{\omega_{KE}^n}.$$

\begin{proof}[Proof of Theorem \ref{volume}]
 Fix a point $q$ such that ${\rm dist}_{\omega_{KE}}(q,\partial\cU)=: 2R(q)$. We construct a cut-off function $\zeta(x) = \eta(\frac{r(x)}{R})  \geq 0$, where $r(x)$ is some smoothening of $d_{\omega_{KE}}(x,q)$ (obtained for instance by running the heat flow) and 
$$ \eta\in [0,1] ,\quad \eta^{-1} (\eta')^2\le C(n), \quad |\eta''| \leq C(n),$$ and let $H = \zeta \varphi$. Note that $\zeta$ satisfies $$\zeta=1~\textnormal{on}~ B_{\omega_{KE}}(q, R),  ~ \zeta=0 ~\textnormal{outside}~B_{\omega_{KE}}(q, 2R).$$ Moreover, by Laplace comparison, we also have that $$\Delta_{\omega_{\KE}}\zeta \geq \frac{-C}{R^2}(1+R).$$ for some dimensional constant $C$. Since both $\omega_{KE}$ and $\omega'_{KE}$ are KE metrics, we have 

$$tr_{\omega_{KE}}\ddbar  \varphi=-n+tr_{\omega_{KE}}\omega'_{KE}\ge n(e^{\frac{\varphi}{n}}-1).$$
Assume $H$ attains a positive maximum at point $Q$ (otherwise $H(q)\le 0$).  Then at point $Q$, we have
\begin{align*}
\Delta H&\ge(\Delta_{\omega_{KE}}\zeta)(\frac{H}{\zeta})+\zeta\Delta_{\omega_{KE} }\varphi+2\Re\langle \nabla\zeta,\nabla\frac{H}{\zeta}\rangle\\
&\ge \frac{-c(n)H}{\zeta}(R^{-2}(1+R))+\zeta n(e^{\frac{\varphi}{n}}-1)+2\Re\langle \nabla\zeta,\frac{1}{\zeta}\nabla H\rangle-2H \langle \nabla\zeta,\frac{1}{\zeta^2}\nabla\zeta\rangle\\
&\ge \frac{-c(n)H}{\zeta}(R^{-2}(1+R))+n\zeta \frac{\varphi^2}{n^2}+2\Re\langle \nabla\zeta,\frac{1}{\zeta}\nabla H\rangle-2 \frac{H}{\zeta R^2}\\
&\ge \frac{-c(n)H}{\zeta}(R^{-2}(1+R)-\frac{H}{c(n)n})+2\Re\langle \nabla\zeta,\frac{1}{\zeta}\nabla H\rangle-2 \frac{H}{\zeta R^2}\\
&=\frac{c(n)H}{\zeta}\Big(-R^{-1}-R^{-2}-\frac{2}{c(n)}R^{-2}+\frac{H}{nc(n)}\Big)+2\Re\langle \nabla\zeta,\frac{1}{\zeta}\nabla H\rangle,
\end{align*}
 where $c(n)$ is a dimentional constant which might change from line to line. By maximum principle on the ball of radius R, noticing that $\nabla H (Q)=0$, we get 
\begin {equation}
H(Q)\le c(n)(\frac{1}{R}+\frac{1}{R^2}).
\end{equation}
Hence $\varphi(q)=H(q)\le 2c(n)\frac{1}{R(q)}$ when $R\geq 1$.  
%Therefore we have 
% $\varphi(q)\le\epsilon$  when ${\rm dist}_{\omega}(q,\partial\cU)\ge\frac{2c(n)}{\epsilon}$.   Switch the role of $\omega$ and $\omega',$ by the same argument, we get $-\epsilon\le \varphi\le +\epsilon$ for $x\in U_{\epsilon}:= \{x|{\rm dist}_{\omega}(x,\partial\cU)\ge\frac{2c(n)}{\epsilon}$ and ${\rm dist}_{\omega'}(x,\partial\cU)\ge\frac{2c(n)}{\epsilon}\}.$  As a conclusion, we prove the theorem.

 \end{proof}
 \begin{remark}
The above theorem is true as long as $\omega_{KE}$ and $\omega'_{KE}$ are complete. No other metric properties of $\omega_{KE},\omega'_{KE}$ are required.
\end{remark}

We prove one corollary of Theorem \ref{volume}.
\begin{corollary}
 Suppose  we are in the setting of Theorem \ref{volume} i.e, $\cU$ admits a complete K\"ahler-Einstein metric $\omega_{KE}$ with negative scalar curvature and $Vol_{\omega_{KE}}(\cU)<\infty$, then for any other complete K\"ahler-Einstein metric $\omega'_{KE}$ with negative scalar curvature, 
 $Vol_{\omega'_{KE}}(\cU)<\infty$.
 \end{corollary}
 \begin{proof}
  it is obvious, since $\varphi=\log\frac{\omega'^n_{KE}}{\omega^n_{KE}}$ is bounded by theorem \ref{volume}.  
 \end{proof}

%Another simple application is a quick proof of uniqueness of complete K\"ahler-Einstein metric on a complex manifold without boundary.
%\begin{corollary} 
%Let $X$ be complex manifold without boundary. Suppose $\omega_{KE},\omega'_{KE}$ are two complete K\"ahler-Einstein metric with negative scalar curvature on $X$, then $\omega_{KE}=c\omega'_{KE}$ for some $c>0$.
%\end{corollary}
%\begin{proof}
% Normalizing the metrics, we may assume that both metrics have constant scalar curvature $-n$. From the K\"ahler-Einstein conditon, we know that
% $$\omega_{KE}=\ddbar\log(\omega_{KE})^n,\omega'_{KE}=\ddbar\log(\omega'_{KE})^n$$
% Hence if we let  $\varphi=\log\frac{\omega'^n_{KE}}{\omega_{KE}^n}$, then
 % $\omega'_{KE}=\omega_{KE}+\ddbar\varphi$ and $\varphi$ satisfies the following equation:
 % $$(\omega_{KE}+\ddbar\varphi)^n=e^{\varphi}\omega_{KE}^n$$
 % From the proof theorem \ref{volume},
 % at a fixed point $p\in X$, as long as we can get a large scale cut off function, which is always true on complete manifold, we have $\varphi(p)<\frac{1}{R}$, for any $R$ large. Hence $\varphi(p)\leq 0$ and similarly we can get $\varphi(p)\geq 0$.  
%\end{proof}

By applying Yau's Schwarz Lemma, we have the following theorem concerning the comparison of two different K\"ahler-Einstein metrics.
\begin{corollary}\label{cmetric}Let $\omega_{KE}$ and $\omega'_{KE}$ be two complete K\"ahler Einstein metrics with negative scalar curvature. 
If moreover the bisectional curvature of $\omega_{KE}$ is smaller than $-K_2.$  where $K_2$ is a positive constant. Then there is a constant $c$ such that $$\frac{1}{C}\omega_\KE'\le \omega_\KE\le C \omega_\KE'.$$ \end{corollary}

\begin{proof}
Let $u=tr_{\omega_\KE'}\omega_\KE$,
by Chern-Lu's inequality we have $$\Delta_{\omega_\KE'} u\ge -K_1u+K_2u^2.$$ where $-K_1$ is the Ricci curvature of $\omega'$. We still use the  cut-off function $\zeta$ as in Theorem \ref{volume}.
Let $G=\zeta u$. Then by the Chern-Lu inequality and the same argument as in the proof of  Theorem \ref{volume}, we get the following inequality
$$G\le \frac{K_1}{K_2}+c(K_1,K_2)R^{-\frac{1}{2}}$$
When $R$ is larger, the estimate is better, hence we have $\omega_\KE\le c\omega_\KE'$ for some constant which depending  on the metric $\omega'$. The volume ratio estimate in Theorem \ref{volume} also gives us the reverse inequality.
\end{proof}

To end this section,  we prove a uniqueness theorem for Dirichlet problem which will be used later in section 4. 
\begin{lemma} \label{unique1} Suppose that there exists a complete K\"ahler metric $\omega=\ddbar \rho$ on $(\cU\setminus p)$. Further assume that the K\"ahler potential   $\rho$ is bounded from above and $\rho(x)\to -\infty$, when $x\to p$.  Then any smooth bounded solution $\varphi$ of the following Dirichlet problem
\begin{equation}\label{unique}
 \begin{cases}
& (\omega+\ddbar\varphi)^n=e^\varphi\omega^n \quad on \quad \cU\setminus p ,\\
&  \varphi|_{\partial\cU}=0.
\end{cases} 
\end{equation}
is unique i.e $\varphi=0$. 
\end{lemma}
 \begin{proof}
  Let $\varphi_{\epsilon}=\varphi-\epsilon\rho$ and $\omega_\epsilon=(1+\epsilon)\omega$, then $\varphi_\epsilon$ satisfies the equation
  \begin{equation}
\begin{cases}
& (\omega_\epsilon+\ddbar\varphi_\epsilon)^n=\frac{e^\varphi}{(1+\epsilon)^n}\omega_\epsilon^n \quad on \quad \cU ,\\
&  \varphi_\epsilon|_{\partial\cU}=-\epsilon\rho.
\end{cases} 
\end{equation}
Since $\rho$ goes to $-\infty$, $\varphi_\epsilon$ admits minimum in $\cU$. If minimum is on the boundary, $\varphi_\epsilon\ge-\epsilon\, \inf_{\partial\cU}\rho$. If the minimum is in the interior point $Q$, 
$\varphi_\epsilon(Q)=\varphi(Q)-\epsilon\rho(Q)\ge\log(1+\epsilon)-\epsilon \max\rho$. In both cases, let $\epsilon\to 0$, we get $\varphi\ge 0$. Similar argument showes that $\varphi\le 0$. Hence $0$ is the unique solution.
\end{proof}
\begin{remark}
 We remark that  Lemma \ref{unique1} can also be proved by using maximum principle when $\omega$ has bounded geometry property.
\end{remark}
\section{K\"ahler-Einstein metrics with bounded geometry} 
\subsection{Preliminaries of the bounded geometry method}
In this section, we want to use bounded geometry methods of \cite{Kob,CY2,TY} to construct complete K\"ahler-Einstein metric on $\cU\setminus \{p\}$. So our geometric domain of interest will be $\cU\setminus \{p\}$ with boundary $\partial\cU$ and also one complete end of infinite distance, which is punctured neighbourhood of $p$ in Euclidean topology.

We recall definitions of quasi-coordinate which are used by \cite{CY,Kob,TY} to deal with complete Riemannian manifolds with bounded curvature but  with shrinking injectivity radius.
\begin{definition} Let $V$ be an open set in $\mathbb C^n$ with coordinates $(v^1,v^2,\cdots,v^n)$.   Let X be an $n$-dimensional
complex manifold and $\phi$ a holomorphic map of $V$ into $X$. $\phi$ is called a
quasi-coordinate map if it is of maximal rank everywhere. In this case,
$(V, \phi,(v^1,v^2,\cdots,v^n))$ is called a quasi-coordinate of $X$.
\end{definition}
\begin{definition}\label{quasi} Let ${\mathcal W}$ be a neighbourhood of $p$ compactly contained in $\cU$, and $\omega$ be a  K\"ahler metric on $(\cU\setminus p)$ which is complete towards $p$. 
A system of quasi-coordinates on $(\widehat{\mathcal W}:={\mathcal W}\backslash\{p\},\omega)$ is a set of quasi-coordinates $\Gamma=\{(V_{\alpha},\phi_{\alpha},(v_{\alpha}^1,v_{\alpha}^2\cdots, v_{\alpha}^n))\}$ of $\widehat{\mathcal W}$ with the following properties:
\begin{enumerate}[label=(\alph*)]
\item $\widehat{\mathcal W}\subset\bigcup_{\alpha}({\rm Image} \,\,{\rm of} \,V_{\alpha})\subset(\cU\setminus p)$;

%\item The complement of certain open neighborhood $U\subset\mathcal U$ of the infinity point $o$ is covered by a finite number of quasi coordinates which are coordinate charts in the usual sense;

\item For each point $x\in \widehat{\mathcal W}$, there is a quasi-coordinate $V_{\beta}$ and $\tilde x\in V_{\beta}$, such that $\phi_{\beta}(\tilde x)=x$ and dist$(\tilde x,\partial V_{\beta})\ge \epsilon_1$ in the euclidean sense, where $\epsilon_1$ is constant independent of $\beta$;

\item There are positive constant $c$ and $A_k,k=1,2,\cdots,$ independent of ${\alpha}$, such that for each quasi coordinate $(V_{\alpha},\phi_{\alpha},(v_{\alpha}^1,v_{\alpha}^2\cdots, v_{\alpha}^n))$,  the following inequalities hold:
$$c^{-1}(\delta_{i\bar j})\le (g_{\alpha i\bar j})\le c(\delta_{i\bar j})$$
$$ |\frac{\partial^{p+q}}{\partial v_{\alpha}^p\bar v_{\alpha}^q}g_{\alpha i\bar j}|<A_{p+q},\forall p,q,$$ 
where $(g_{\alpha i\bar j})$ denote the metric tensor of the Riemannian metric associated to $\omega$ on chart $(V_{\alpha},\phi_{\alpha},(v_{\alpha}^1,v_{\alpha}^2\cdots, v_{\alpha}^n))$.
\end{enumerate}
\end{definition}

Roughly speaking, a set of quasi coordinates of metric domain $(\widehat{\mathcal W}, \omega)$ is a set of coverings of $\mathcal W$ by charts such that the pull back of $\omega$ satisfies uniform bounded metric properties.
Now we define the Cheng-Yau function space.
\begin{definition}\label{C-Y}
We define the H\"older space of $C^{k,\alpha}$ function on $\widehat{\cU}:=\cU\setminus p$ by exploiting the quasi-coordinate system. For any nonnegative integer $k,\alpha\in (0,1)$,  we define 
$$\|u\|_{k,\alpha}(\widehat{\mathcal W})=\sup_{V_{\beta}\in \Gamma } (\sup_{z\in V_{\beta}}\sum_{p+q\le k}|\frac{\partial^{p+q}}{\partial v_{\beta}^p\partial \bar v_{\beta}^q}u(z)|+\sup_{z,z'\in V_{\beta}}\sum_{p+q=k}|z-z'|^{-\alpha}|\frac{\partial^{p+q}}{\partial v_{\beta}^p\partial \bar v_{\beta}^q}u(z)-\frac{\partial^{p+q}}{\partial v_{\beta}^p\partial \bar v_{\beta}^q}u(z')|)$$
let us introduce one more compact set $V$ with $\widehat{\cU}\setminus\widehat{\mathcal W}\subset V\subset\widehat{\cU}$ to cover whole $\widehat{\cU}$. Now define:
$$\|u\|_{k,\alpha}(\widehat{\cU})=\|u\|_{k,\alpha}(\widehat{\mathcal W})+\|u\|_{k,\alpha}(V)$$
\end{definition}
The function space $C^{k,\alpha}(\widehat{\cU})$ is the completion of  $\{u\in C^k(\widehat{\cU})| \|u\|_{k,\alpha}(\widehat{\cU})< \infty\}.$

 \begin{remark}
 The existence of quasi coordinate is crucially used in our proof.  The classical interior Schauder estimate for a linear elliptic operator $L,$ is as follows:
$$\|u\|_{C^{k,\alpha}(V_1)}\le C(\sup|u|_{V_2}+\|Lu\|_{C^{k-2,\alpha}(V_2)}),$$ where $V_1\subset\subset V_2\subset R^m$.
 Notice that the constant $C$ depends on the ellipticity of $L,$ the $C^{k-2,\alpha}$ norms of the coefficients of $L$ and the distance between  $V_1$ and $\partial V_2$. If we have a quasi coordinate system defined above, the Schauder estimate on $\widehat{\cU}$ is reduced to that on a fixed bounded domain in Euclidean space. 
\end{remark}
We also introduce the following assumption  for our geometric domain ($\cU,p$).
\begin{definition}
{\bf Bounded geometry of order $k,$}\label{property A}: Let ($X,p$) a germ of isolated log canonical singularity embedded in ($\mathbb C^N,0$), we call ($X, p$)
has bounded geometry of order $k$. If
\begin{enumerate}\item  there is a  complete metric $\omega=\ddbar\rho$ defined on ($X\setminus p$) which  has a system of quasi coordinates up to $k$th derivative of metric $g$, see definition \eqref{quasi} item $(d)$. 
\item  there is a function $M$ on $(\cU\setminus p)$ satisfying $Ric(\omega)+\omega=\ddbar M$ and for any $i\le k, ||\na_{\omega}^i M||<C_i.$ (Here the potential function $M$ is not unique, we only require one of them satisfy the boundedness property, and in this note, the most interesting case is $M=0$). 
\end{enumerate}
\end{definition}

Before we proceed, we state and prove  the following modified version of Yau's generalized maximum principle on noncompact manifold.
\begin{lemma}\cite{CY}\label{max}
Suppose $(\widehat{\cU},\omega)$ is of bounded geometry of order $k$ with $k\geq 2$. Let $ f$ be a smooth function on $\widehat{\cU},$ which is bounded from above, and $\sup f> \sup_{\partial{\mathcal U}}f,$ then there is a sequence $\{y_i\}$ in $\widehat{\cU}$ such that $\lim_{i\rightarrow\infty} f(y_i)=\sup f, \lim_{i\rightarrow\infty} |\na_g f|(y_i)=0$
and $\lim_{i\rightarrow \infty} |\Delta_gf|(y_i)\leq 0,$ where the derivatives  are taken with respect to metric $g$ associated to the the K\"aher form $\omega.$
\end{lemma}
\begin{proof}: Without loss of generality, let us assume that $\sup f=0$.
If  $\sup f$ is attained, the lemma is obvious.   Otherwise we choose a sequence $x_i$ with $\lim f(x_i)=0.$ It is easy to see $\{x_i\}$ must go to infinity. Now at each point we take a quasi coordinate chart $V_i$ covering $x_i.$ On each $V_i$, define a non-negative  function $\beta^i: V_i\rightarrow \mathbb R$ such that $$\beta^i(x_i)=1,\beta^i=0 \text{ on }\partial V_i,~ \beta^i\le C,|\nabla\beta^i|\le C, \text{ and }(\beta^i_{p\bar q})\ge -C(\delta_{p\bar q}),$$
where $C$ is positive number independent of $i,$ and all norms are taken with respect to the Euclidean norm. 
Now consider $$\frac{-f}{\beta^i}$$ as a function on $V_i.$ Notice that $\frac{-f}{\beta^i}$ blows up on the boundary of $V_i$, so it admits minimum at  point $y_i$ which is in the interior of $V_i$. Now let 
$$\frac{-f}{\beta^i}(y_i)=\inf_{V_i}\frac{-f}{\beta^i}.$$
Then
\begin{align*}
\frac{-f}{\beta^i}(y_i)&\le\frac{-f}{\beta^i}(x_i)=-f(x_i),\\
\frac{df}{f}(y_i)&=\frac{d\beta^i}{\beta^i}(y_i),\\
\frac{f_{p\bar q}}{f}(y_i)&\ge\frac{\beta^i_{p\bar q}}{\beta^i}(y_i).
\end{align*}
Using these inequalities and our choice  of $\beta^i$, we have
\begin{align*}
0<-f(y_i)&\le -Cf(x_i),\\
|df(y_i)|&\le -Cf(x_i),\\
(f_{p\bar q})(y_i)&\le -Cf(x_i)(\delta_{p\bar q}).
\end{align*}
By the bounded geometry of quasi coordinates, the above norms can also be take with respect to the metric $\omega.$ Hence sequence $\{y_i\}$ satisfies all the properties required in the lemma.
\end{proof}
\subsection{Construction of K\"ahler-Einstein metrics with bounded geometry} Using function $\rho$ in item (1) of Definition (\ref{property A}) to define a domain:  $$(\cU\setminus p):=\{\rho<a\}$$ 
We point out that function $\rho$ used here is different from the one used in definition (\ref{domain1}). The main goal of this section is to prove the following theorem concerning the solvability of K\"ahler-Einstein equation on $(\cU\setminus p)$ by a perturbative method.
\begin{theorem}\label{Dirichelt-bounededgeometry}Suppose $(X,p)$ is a germ of log canonical singularity and a punctured neighbourhood $\cU\setminus p$ of $p$ admits a complete K\"ahler metric $\omega$ with bounded geometry of order $k$. 
Then for any smooth function $\psi$ on the boundary $\partial\cU$, the following Dirichlet problem 
\begin{equation}\label{Dirichlet}
\begin{cases}
& ( \omega + \ddbar \varphi ) ^n =e^{\varphi+M} {\omega^n} \,\,\,{\rm on }\, \,\,\cU, \\
&  \varphi|_{\partial\cU}=\psi.
\end{cases}
\end{equation}
admits a solution in Cheng-Yau function space defined in \eqref{C-Y} with $||\varphi||_{k,\alpha}<C(\psi,\rho,M)$. 
\end{theorem}
We first take the function space $U$ to be an open set of $C^{k,\alpha}(\cU)$, which is defined in \eqref{C-Y}, as follows:
$$U=\{\phi\in C^{k,\alpha}(\cU): \frac{1}{c}(g_{\alpha i\bar j})\le (g_{\alpha i\bar j}+\phi_{i\bar j})\le c(g_{\alpha i\bar j}), {\rm in\,\, each\,\, quasi\,\, coordinate\, } V_{\alpha }\},$$ 
for some  constant $c$, which however is not fixed. 

\begin{proof}[Proof of theorem \ref{Dirichelt-bounededgeometry}:] The proof consists of several steps.\\

\item{\bf Step 1: Find an $\omega$-PSH extension of function $\psi$ to the domain $\cU$. } Choose an arbitrary smooth extension $\psi_1$ of $\psi$, which is supported on  $\{x:a\ge\rho(x)\ge c\}$ Also choose a convex monotone increasing  function $H:[-\infty,a]\rightarrow \mathbb R$ which is zero on $[-\infty,b]$ for some constant $b<c$.
Now define  $$P=AH(\rho)-AH(a).$$ By chossing A  large,   $\omega+\ddbar P+\ddbar\psi_1$  is still a K\"ahler form. Choosing  $P+\psi_1$ as our new extension of $\psi$, then by  construction  $\ddbar (P+\psi_1)$ is supported on a neighborhood of $\partial\cU $. Therefore $\omega$ coincides with $\omega+\ddbar P+\ddbar\psi_1$ in a punctured neighbourhood of $p$.  

In sum,  we have
\begin{equation}
\begin{cases}
& ( \omega+ \ddbar (\psi_1+P) ) ^n =e^{-F} {\omega^n}  \,\,\,{\rm on }\, \,\,\cU\\
& \psi_1+P|_{\partial\cU}=\psi .
\end{cases}
\end{equation}
where the function $F$ is still in Cheng-Yau's function space. Hence if we define  $\tilde\omega$ by $\tilde\omega=\omega+\ddbar(\psi_1+P),\tilde\varphi=\varphi-(\psi_1+P)$ and $\tilde F=F+\psi_1+P.$
Simple calculation shows the Equation \eqref{Dirichlet} is equivalent to 
\begin{equation}
\begin{cases}
& ( \tilde\omega  +\ddbar\tilde\varphi) ^n =e^{\tilde\varphi+\tilde F} {\tilde \omega}^n \,\,\,{\rm on }\, \,\,\cU \\
& \tilde\varphi|_{\partial\cU}=0 .
\end{cases} 
\end{equation}

So from now on, we will focus on zero boundary value problem.

The rest of the proof  is by continuity method, which is based on a combination of  estimates from \cite {CKNS,Blo,CY,Kob}.
We set up the continuity method as follows:
\begin{equation}\label{Con}
\begin{cases}
& ( \omega + \ddbar \varphi_{t} ) ^n =e^{\varphi_t+tM} {\omega^n} \\
& \varphi_t|_{\partial\cU}=0.
\end{cases} 
\end{equation}
where $M$  belongs to Cheng-Yau function space defined in \eqref{C-Y}.\\

\item{\bf Step 2: Openness part in the continuity method.}
As usual, the openness will follow from the inverse mapping theorem. We need to show the linearized equation at $\omega_t$
\begin{equation}\label{Lin}
\begin{cases}
& \Delta_{\omega_t} h-h=v\,\,\,{\rm on }\,\,\, \cU\\
& h|_{\partial\cU}=0,
\end{cases}
\end{equation}
has a unique solution in $C^{k,\alpha}(\cU)$ with the estimate $$\|h\|_{k,\alpha}(\cU)\le c\|v\|_{k-2,\alpha}(\cU)$$
for some constant $c$ independent of the function $v$.\\
We first remark here that $\omega_t:=\omega+\ddbar\varphi_t$ is a complete metric of bounded geometry up to $k-2$ covariant derivatives by  the function choice of function space $U$ at the beginning of the proof.  
Next take an exhaustion $\{\cU_i\}$ of the domain $\cU$ towards the infinity. (Here the boundary of our compact domain $\cU_i$ has two components and one of them coincide with $\partial \cU$).  
Following equation
\begin{equation}
\begin{cases}
& \Delta_{\omega_t} h_i-h_i=v \,\,\, on \,\,\, \cU_i \\
& h|_{\partial\cU_i}=0,
\end{cases} 
\end{equation}
has a unique solution $h_i$ for each $i$. Maximum principle implies that $\sup_{\cU_i}|h_i|\le  \sup |v|.$
Interior Schauder estimate of our function space implies that 
$$\|h_i\|_{k,\alpha}(K)\le c\|v\|_{k-2,\alpha}(\cU).$$ 
for any compact set $K$ strict away from $\partial\cU$. This inequality, combined with standard global Schauder estimate for a fixed compact set $V$ containing  $\partial\cU$ imply that  $h_i\rightarrow h$ pointwise with $$\Delta_{\omega_t}h-h=v, h|_{\partial\cU}=0$$
Moreover, we have
$$\|h\|_{k,\alpha}(\cU)\le c\|v\|_{k-2,\alpha}(\cU).$$  Hence we establish the openness part.\\

\item{\bf Step 3: $C^0 $estimate.}
We have the following equality:
\begin{align*}
\varphi +M&=\log {\det(g_{i\bar j}+\varphi_{i\bar j})}-\log {\det g_{i\bar j}}\\&=\int_0^1\frac{\partial}{\partial t} \log{\det(g_{i\bar j}+t\varphi_{i\bar j})}dt\\&=\int^1_0B(t)_{i\bar j}\varphi_{i\bar j}dt, 
\end{align*}
where $B(t)_{i\bar j}$ is the cofactor matrix of matrix $(g_{i\bar j}+t\varphi_{i\bar j})$. Since at a point $x\in \Omega$ we may asume $g_{i\bar j}=\delta_{i \bar j}$ and $\varphi_{i\bar j}=\delta_{i\bar j}\varphi_{i\bar i}$. If $\varphi_{i\bar i}\geq 0$, then $\frac{\varphi_{i\bar i}}{1+t\varphi_{i\bar i}}\leq \varphi_{i\bar i}$. If $\varphi_{i\bar i}\leq 0$, then $\frac{\varphi_{i\bar i}}{1+t\varphi_{i\bar i}}\leq \varphi_{i\bar i}$.  Hence we have two inequalities as follows:
$$\varphi+M\le\Delta_{\omega}\varphi,$$
$$\varphi+M\ge\Delta_{\omega_1}\varphi.$$
where $\omega_1=\omega+\ddbar \varphi$. By  Lemma \ref{max}, we get the $C^0$ estimate.\\
%We remark here in our situation, we have boundary and Yau is considering a complete manifold ,but we can easily verify that if $\sup\varphi>\varphi|\Omega,$
%then the sequence of points chosen by yau should go away from the boundary.\\

\item {\bf Step 4: $C^1$ boundary estimate.} On the one hand since  
$\varphi+M\le\Delta_{\omega}\varphi$,  we construct a barrier function $h$ from above as follows. Take a domain $V\subset \cU$ satisfying $p\notin V$. We also require that  $\partial V=\partial \cU\cup C$, where $C$ is a  smooth connected manifold disjoint with $\partial\cU$.  Then derive $h$  by solving the following Dirichlet problem in $V$:
\begin{equation}
\begin{cases}
& \Delta_{\omega} h=c,  \\
& h|_{\partial \cU}=0 \,\,{\rm and}\,\,h|_{C}=d.
\end{cases}
\end{equation}
where $d:= \sup_{\cU} |\varphi|$ and $c:=\inf_{\cU} (\varphi+M).$
Then maximum principle implies that $h\ge\varphi$ in $V$.  
%Then the normal derivative of $h-\varphi$ on the boundary $\partial\cU$ is nonpositive which gives a lower bound for the normal derivative of $\varphi$ on $\partial\cU$.

On the other hand, we construct a barrier function $h_1$ from below as follows. Take the global strict $\omega$  -PSH function $P$ we constructed in step 1 and choose a constant $B$ large enough such that
\begin{equation}
\begin{cases}
&(\omega+\ddbar BP)^n\ge e^{sup\varphi+M}\omega^n \,\,{\rm on}\,\,\,\,V,  \\
& BP\le\varphi \,\,{\rm on}\,\, C \,\,{and } \,\,BP=0\,\,{\rm on}\,\,\partial\cU.\\
\end{cases} 
\end{equation}
Then maximum principle of Monge-Amp\`ere equation implies  $h_1:=BP\le\varphi$ on $V$. 
%Hence the normal derivative of $h-\varphi$ on the boundary $\partial\cU$ is nonnegative which gives a upper bound for the normal derivative of $\varphi.$
Noticing that $h$ and $h_1$ coincide with $\varphi$ on $\partial \Omega$, we get the boundary gradient estimate of $\varphi.$\\

\item {\bf Step 5: Global $C^1$ estimate.}  Since  on noncompact manifold, we do not necessary have maximum point with gradient vanish etc. We will follow  \cite{Blo} to obtain the  $C^1$ estimate. 
Define $\phi=\log|\nabla\varphi|^2-\gamma(\varphi)$ where  $\gamma$ is monotone increasing function to be determined.
Assume that $\sup_{\Omega} \phi$ is not obtained on $\partial\cU$, then by the generalized maximum principle, we can find a point $q\in\Omega$ with $$\phi(q)+\epsilon>\sup_{\Omega} \phi, |V|_{\omega}(q)<\epsilon,\Delta_{\omega}\phi(q)<\epsilon$$
where $V:=\nabla_{\omega}\phi.$ By \eqref{PS}, at $q$, we have
\begin{equation*}
\Delta'\log \,|\nabla\varphi|_g^2
\ge
\frac{2 {\rm Re}\nabla_m\log F \nabla^m \varphi
}{ |\nabla\varphi|_g^2}
-\Lambda \,tr_{g'}g
+
2\Re\langle \frac{\na|\na\varphi|_g^2}{ |\na\varphi|_g^2}, \frac{\na \varphi}{ |\na\varphi|_g^2}\rangle_{g'} - 
2\Re\langle \frac{\na|\na\varphi|_g^2}{ |\na\varphi|_g^2}, \frac{\na \varphi}{ |\na\varphi|_g^2}\rangle_g 
\end{equation*}
Now using $V=\nabla\log|\nabla\varphi|^2-\gamma'\nabla\varphi$, we get 
\begin{align*}
2\Re\langle \frac{\na|\na\varphi|_g^2}{ |\na\varphi|_g^2}, \frac{\na \varphi}{ |\na\varphi|_g^2}\rangle_{g'} - 
2\Re\langle \frac{\na|\na\varphi|_g^2}{ |\na\varphi|_g^2}, \frac{\na \varphi}{ |\na\varphi|_g^2}\rangle_g
&=2\Re\langle {V+\gamma'\na\varphi}, \frac{\na \varphi}{ |\na\varphi|_g^2}\rangle_{g'} - 
2\Re\langle V+\gamma'\na\varphi, \frac{\na \varphi}{ |\na\varphi|_g^2}\rangle_g\\&>-\epsilon tr_{g'}g-\epsilon-\gamma'
\end{align*}
The inequality above with 
$$\Delta'\gamma(\varphi)=\gamma'\Delta'\varphi+(\gamma'')|\na\varphi|^2$$
imply that
$$\epsilon tr_{g'}g>\Delta'\phi>(\gamma'-\epsilon-\Lambda)tr_{g'}g+(\gamma'')|\na\varphi|^2-\gamma'-\epsilon-n-C$$
where $C$ is the bound of gradient of function $\log F$. Now we construct our function $\gamma$ as
$$\gamma(x)=(\Lambda+2)x-\frac{1}{x+C'+1}$$
where $C'$ is the lower bound of $\varphi$. Then by standard argument we get global $C^1$ estimate.\\

\item {\bf Step 6: Boundary $C^2$ estimate.} We notice the argument of \cite {CKNS} of  boundary $C^2$ estimate is purely local around the boundary and  our equation can be written as

\begin{equation}
\begin{cases}
& \det\varphi_{i\bar j}=e^{\varphi+f} \,\,\, on \,\,\, \Omega \\
&  \varphi|_{\partial\cU}=\psi,
\end{cases} 
\end{equation}
locally, this is exactly one of the equation considered in \cite{CKNS}, hence the estimate follows by the fact that $\ddbar\rho$ is strictly positive in a neighbourhood of $\partial\cU$.\\

\item {\bf Step 7: Global  $C^2$ estimate and completion of the proof.}
We have the well-known inequality as follows:
$$\Delta'\log tr_g g'\ge-Btr_{g'}g-C$$ where B,C depends on the geometry of good background metric $g$ and Ricci curvature of the volume form on the right hand side of  the equation.
Notice that $$\Delta'\varphi=n-tr_{g'}g.$$
By setting A=B+C+1,we have the differential inequality $$\Delta'(\log tr_gg'-A\varphi)=tr_{g'}g-An.$$ 
This inequality and  the boundary $C^2$ estimate imply the global $C^2$ estimate. Then by 
Evans-Krylov Theorem, we have interior $C^{2,\alpha}$ estimate. Now the metric is uniformly bounded in a neighbourhood of the boundary,  the local argument of \cite{CKNS} will give global $C^{2,\alpha}$ estimate.

\end{proof}
\begin{remark}
The $C^{2,\alpha}$ estimate of $\varphi$ depends on the extension  $\psi_1$ of  our boundary function  $\psi$ in Step 1, it is not clear if we can have $\|\psi_1\|_{2,\alpha}<C\|\psi\|_{2,\alpha}$ for some constant $C$ independent of the boundary value.
\end{remark}

%\begin{remark}
%We use the decay of $\varphi$ to prove any complete K\"ahler-Einstein $\omega'$ comes from the solution of (\ref{Dirichlet}). Even if we first assume $\omega'$ comes from one solution of equation (\ref{Dirichlet}) corresponding to a choice of boundary function $\psi$, since we can have very huge perturbation of boundary condition $\psi$, from which we can only conclude boundedness of $\varphi$ globally on $(\Omega\setminus p)$. Hence our proof of volume rigidity is more like a geometric argument rather a PDE argument.
%\end{remark}

In the following lemma, we prove that any complete local K\"ahler-Einstein metric is indeed obtained from solving Monge-Amp\`ere equations (\ref{Dirichelt-bounededgeometry}) with different Dirichlet boundary conditions. 
\begin{lemma}\label{ddbar}
Let $\omega_{KE}$ be a complete K\"ahler-Einstein metric with bounded geometry property on $\cU\setminus p$ and $\omega'_{KE}$ be another complete K\"aher-Einstein metric on $\cU\setminus p$. Then there exists a function $\psi$ on $\partial\cU$ and a solution $\varphi$ of equation (\ref{Dirichlet}) (with $\omega=\omega_{KE},M=0$,Dirichlet boundary value $\psi$) obtained by the bounded geometry method as in Theorem \ref{Dirichelt-bounededgeometry}. Moreover, $\omega'_{KE}=\omega_{KE}+\ddbar\varphi$.
\end{lemma}
\begin{proof}
 From the fact that $\omega'_{KE}$ is a complete K\"ahler-Einstein metric and Theorem \ref{volume}, we know that $\omega'_{KE}=\omega_{KE}+\ddbar\tilde\varphi$ and $\tilde\varphi$ is a bounded smooth  solution of equation (\ref{Dirichlet}). On the other hand, by Theorem \ref{Dirichlet}, we can find another function $\varphi$ solving equation (\ref{Dirichlet}) 
 %$\tilde\omega'_{KE}=\omega_{KE}+\ddbar\tilde\varphi$
with  Dirichlet boundary   value $\psi:=\tilde\varphi|_{\partial\cU}$.
By the uniqueness Lemma \ref{unique1}, we conclude that $\varphi=\tilde\varphi,$ and hence $\omega'_{KE}=\omega_{KE}+\ddbar\varphi$. This completes the proof.
\end{proof}

We shall now prove Theorem \ref{main2} on the asymptotic behaviour of KE metrics we constructed in Theorem \ref{Dirichelt-bounededgeometry}. 
With the help of estimates of higher order derivatives of $\varphi$ and quasi coordinates, we prove the stability of local complete K\"ahler-Einstein metrics with negative scalar curvature near isolated log canonical singularity.
\begin{proof}[Proof of Theorem \ref{main2}]:
 First of all, since $\omega'_{KE}$ is complete, by Lemma \ref{ddbar}, $\omega_{KE}'=\omega_{KE}+\ddbar\varphi$ where $\varphi$ is the solution from Theorem \ref{Dirichelt-bounededgeometry} with suitable Dirichlet boundary value. For any point $q\in (\cU\setminus p)$, we can  choose a quasi coordinate $(\widehat V,\phi)$ covering $q$ such that there is a point  $\widehat q\in V\subset \widehat V$, $\phi(\widehat q)=q$ and ${\rm dist}(\widehat q,\partial \widehat V)\ge {\rm dist}(V,\partial \widehat V)\ge{\epsilon_1}.$ Let $\beta$ be the cut-off function we constructed in the proof of Lemma \ref{max}. Then we have the following inequalities:
\begin{center}
$\sum_{i=1}^{k}|\beta^{(k)}|_{\rm Euc}\le B_k,$
\end{center}
where $B_k's$ are  constants independent of $p$ and $V$ by the existence of quasi coordinates \ref{quasi}.  This is true because under the construction of the system of quasi coordinates, 
%the source domains are all identical to $\widehat\Omega$ and the point we are considering is alway sitting inside a fixed subdomain $\Omega\subset\widehat\Omega$ with 
we have ${\rm dist}(V, \partial\widehat V)\ge\epsilon_1>0$, hence controlling the derivatives of cut-off function uniformly. Actually, we can even assume the covering domains we choose are $B_{\frac{1}{4}\epsilon_1},B_{\frac{1}{2}\epsilon_1}$ by subdividing the original coverings.
By our previous proof of the a priori estimates of $\varphi$ from equation \ref{Dirichlet}, we also have the following inequalities, for any point $q\in (\cU\setminus p)$ and any nonnegative integer $k$:
\begin{center}\label{CK}
$\sum_{i=1}^{k}\|\nabla^{(k)} \varphi\|_{\omega_{KE}}(q)\le C_k.$
\end{center}

 When $k=0$, the $C^0$ decay of $\varphi$  is  Theorem \ref{volume}.  For $k\ge 1,$
 we do computations in the quasi coordinate as follows:
$$-\int_{\widehat V}\beta\varphi\Delta\varphi=\int_{\widehat V}\beta|\nabla\varphi|^2+\int_{\widehat V}\varphi\langle\nabla\varphi,\nabla\beta\rangle,$$
$$\int_{V}|\nabla\varphi|^2\le\int_{\widehat V}\beta|\nabla\varphi|^2\le C\cdot (C_2+ B_1)\cdot \sup_{\widehat V}|\varphi|.$$
Similarly, using integration by part,
$$\int_V \beta|\nabla^{k}\varphi|^2=\int_V(\sum_{i+j=2k-1}\nabla^i\beta\nabla^j\varphi)\varphi\le  (B_k+ C_{2k-1})\cdot \sup_{\widehat V}|\varphi|.$$
 Note that by theorem \ref{Dirichelt-bounededgeometry}, there is constant $C(\psi,\omega_{KE},n) $ such that $\frac{\omega_{KE}}{C(\psi,\omega,n)}\le \omega'_{KE}\le C(\psi,\omega_{KE},n)\omega_{KE}$. Suppose ${\rm dist}_{\omega_{KE}}(q,\partial\cU)\ge R$, then ${\rm dist}_{\omega'_{KE}}(q,\partial\cU)\ge \frac{1}{C(\omega_{KE},\psi,n)}R,$
by the triangle inequality, we have $$dist_{\omega_{KE}}(\partial(\phi(\hat V)),\partial\cU)\ge dist_{\omega_{KE}}(q,\partial\cU)-dist_{\omega_{KE}}(q,\partial(\phi(\hat V)))\ge R-C'\epsilon_1\ge\frac{R}{2}.$$
where $\widehat V$ is a covering of point $q$ and $C'$ is the metric equivalence constant in the definition of quasi coordinates \eqref{quasi}, where depends on the geometry of $\omega_{KE}$. Similarly, we have  $$dist_{\omega'_{KE}}(\partial(\phi(\hat V)),\partial\cU)\ge\frac{R}{2C(\omega_{KE},\psi,n)}.$$ Hence by the $C^0$ estimate of $\varphi$ in Theorem \ref{volume},
$$\int_{V}|\nabla^{k}\varphi|^2\le C(\psi,\omega_{KE},n)\cdot (C_{2k-1}+B_{k})\cdot\frac{c(n)}{R}.$$
Now that we have $L^2$ norm control of all higher order  derivatives, by Sobolev embedding on Euclidean space and property $(d)$ of quasi-coordinate,  we can conclude that $\sum_{i=1}^{k}\|\nabla^i\varphi\|_{\omega_{KE}}(q)\le \frac{1}{R(q)^{1/2}}\cdot C(k,\omega_{KE}, \psi)$ 
\end{proof}
\begin{remark}We remark that the decay rate $R^{\frac{-1}{2}}$ of function $\varphi$ obtained in the above theorem is far from optimal. The optimal decay rate might need a case by case treatment depending on the type of the singularity.  It is pointed out to us by Professor Hein that if $|\varphi|\leq R^{-1}$,  and if all  its higher order derivatives are bounded, then by using quasi coordinates and the regularity theory of Monge-Amp\`ere equations, the decay rate can be improved to $R^{-1}$. 
\end{remark}
We use the metric stability result to show that to construct a geometric domain with bounded geometry property in a punctured neighbourhood of an isolated log canonical singularity, we only need to show there exists a K\"ahler metric with bounded geometry property on its finite cover. Therefore, once for some singularity has bounded geometry property, so it does for its finite quotient, hence more example of singularity with bounded geometry property  will be obtained. This will be useful when the metric on the finite cover is not necessary invariant under the finite group action.  
\begin{corollary}
 \label{quotient}Let $(X,p)$ be an isolated log canonical singularity embedded in $\mathbb C^N$ and $(\cU\setminus p, \omega=\ddbar\rho)$ is K\"ahler-Einstein  domain with bounded geometry property of infinite order defined  as in Definition \ref{property A}.  Let $(Y,p')$ be the quotient space  $(X,p)/ G$, where $G$ is a finite group acting freely on $(X\setminus p)$ and fixing point $p$.  Then for any $k\geq 1$, there is a punctured neighbourhood $\cU'\setminus p'$ of $p'$ admitting a K\"ahler metric $\omega'$ with bounded geometry property of order $k$. Hence Theorem \ref{main2} can be applied to the finite quotient domain $(\cU'\setminus p',\omega')$.   
\end{corollary}
%\begin{remark}We remark that for one type of  isolated log canonical singularity in lemma \ref{quasi1}, i.e the $\mathbb Z_2$ quotient of cusp singularity, (see \cite{Kob2} page 333 for details), due to the limited knowledge of the authors, the model metric  is only constructed on the covering and we are unable to find the  $\mathbb Z_2$ action, which makes the metric invariant under the group action.  Hence we are unable to construct a good model metric by hand explicitly.   But we may use the corollary above to cover this type of surface log canonical singularity. 
%\end{remark}
\begin{proof}
For simplicity, assume $G=\mathbb Z_2$. Let $f$ be the non trivial element of $G$. Define
$$\rho':=\rho+f^*\rho,\omega':=\ddbar\rho'$$
then $\rho'$ is invariant under the $\mathbb Z_2$ action. Note that $\hat\omega=\ddbar f^*\rho$ is also a complete K\"ahler-Einstein metric with bounded geometry property near $p$, hence by Theorem \ref{main2}, we have
$$\hat\omega=\omega+\ddbar\varphi,\sum_{i=1}^{k}\|\nabla^i\varphi\|_{\omega}(q)\le  \frac{C(k,\omega, \psi)}{R^{1/2}(q)}.$$ In particular, this shows that $\omega'=\omega+\frac{1}{2}\ddbar\varphi$ has bounded geometry property of order $k$  by possibly shrinking the domain. This completes the proof. 
\end{proof}

\section{Log canonical singularities on surfaces}
\subsection{Model metrics with bounded geometry property}
In this short subsection, we provide some explicit examples of $(X,p)$ with bounded geometry property of infinite order. We mainly focus on complex dimension 2 and expect that there are  more examples in higher dimension by using arithmetic quotient of  symmetric domains (cf. \cite{Fang-Fu,YZ}). We should remark that all log canonical surface singularity admits a model metric with bounded geometry property. However, in  complex dimension higher than two, there are  examples of log canonical singularity which  admit a local K\"ahler-Einstein metrics but these metrics do not have bounded geometry property. For such an  example,  we can take a negative line bundle over a non flat Calabi-Yau manifold and then contract the zero section. 
% Uniformizatoin of  2 dimensional isolated log canonical singularity by K\"ahler-Einstein metric.

\begin{lemma}\cite{Kob2, Nak, Ye}\label{quasi1}\\
Any isolated normal log canonical (non log terminal) surface singularity can be uniformized by bounded symmetric domains with invariant K\"ahler-Einstein metric $\omega=\ddbar\rho$ and classified as:
\begin{enumerate}
%\item[(1)]$\mathbb C^2/\Gamma,\Gamma$ a finite subgroup of $\mathbb U(2)$ possibly containing reflections. Invariant metric is defined by 
%$\ddbar \rho,\rho=\log(\frac{1}{(1-\|Z\|^2)^2})$

\item[(1)]One point partial compactification of $\mathbb B^2/\Gamma,\Gamma$ a parabolic discrete subgroup of $Aut(\mathbb B^2)$.
Invariant metric is defined by $\ddbar \rho,\rho=\log\frac{1}{1-|z_1|^2-|z_2|^2}$, where $z_1,z_2$ are complex coordinates of $\mathbb B^2$
% \textcolor{red}{Wrong here}.
\item[(2)]One point compactification of $\mathbb H \times \mathbb H/\Gamma,\Gamma$ a parabolic discrete subgroup of $Aut(\mathbb H\times \mathbb H)$ corresponding to a boundary point.
Invariant metric is defined by $\ddbar \rho,\rho=\log(\frac{1}{y_1y_2})$, where $y_1,y_2$ are the imaginary part of the complex coordinates of $\mathbb H\times \mathbb H$.
 \end{enumerate}

\end{lemma}
\begin{remark}
Type 1 singularities, up to a finite quotient, have a minimal resolution whose exceptional divisor is an elliptic curve. Type 2 singularities, once again up to finite quotients, have  minimal resolutions whose exceptional divisor  is a cycle of smooth rational curves.
 %We include the nice construction in the appendix for reader's convenience.
\end{remark}
%\begin{remark}
%Another interesting example of isolated log canonical singularity uniformized by the Bergman metric is proved in \cite{Fang-Fu}. 

%Let  $A$ be an abelian variety with complex dimension $n$ and $N$ a negative line bundle on $A.$ By contracting the zero section of $N,$  one obtains a singular variety $\widehat X$. 
%with the contraction map $f:\widehat M\rightarrow \widehat X$.  % restricts to a natural resolution of singularities of $(\widehat X,o)$, where $(\widehat X,o)$ denotes the germ of the isolated singularity. 
%Let   $(\widehat X,o)$ denote the germ of the isolated singularity of $\widehat X$.  Then there is an open neighborhood (in Euclidean topology) $\widehat {\mathcal U}$ of $o$ in $\widehat X$ such that ${\mathcal U}:=\widehat{\mathcal U}\backslash \{o\}$  is a smooth quotient space of a unit complex ball $\mathbb B^{n+1}$ by a discrete subgroup of ${\rm Aut}(\mathbb B^{n+1}).$ As a consequence, $\mathcal U$  has a negative K\"ahler- Einstein metric  induced from the Bergman metric of the ball which is complete towards $o$.

%\end{remark}

\begin{remark} The invariant K\"ahler-Einstein metrics in Lemma \ref{quasi1} have a system of quasi coordinates in a punctured neighborhood of the isolated  log canonical singularities. This is the main property we will use in the following proof of Corollary \ref{surface}.
\end{remark}

\begin{proof}[Proof of Corollary \ref{surface}] By Lemma \ref{quasi},
for any isolated log canonical surface singularity, there is a reference metric with bounded geometry of order $k$  for any $k$. Hence Corollary \ref{surface} is a direct consequence of theorem \ref{main2} except for the statement about the Gromov-Hausdorff convergence. We deal with each of the cases separately. 
\begin{enumerate}
\item[\textnormal{(i)}] By \cite{Fang-Fu}, the uniformization metric in the upper half plane model is $$\omega=-\ddbar(\textnormal{Im} u-|v|^2)=-\ddbar\log(-\log|\sigma_D|_{h_D}),$$ where $D$ is an elliptic curve. The total space of resolution is a negative line bundle over the elliptic curve.  Direct calculation shows that the metric $\omega$ degenerates along the tangential direction of $D$ and the $S^1$ circle direction contained in the fiber of the negative line bundle. Hence it is clear that, when we choose the base point $p_j\to p$, $(\cU, p_j, g_{KE})\to \mathbb R$ in the pointed Gromov-Hausdorff topology.
\item[\textnormal{(ii)}] 
In this case,  the model is $\mathbb H\times \mathbb H/\Gamma$ with the invariant metric given by $\ddbar\log(\frac{1}{y_1y_2})$ where $y_i$ is the imaginary part $z_i,i=1,2$. To make things clear, we describe what is the group action (see \cite{Kob} page 55 or \cite{Kob2} page 344 for more details).

 Let $G(M,V)=\bigl\{\bigl(\begin{smallmatrix}\epsilon&\mu\\0&1\end{smallmatrix}\bigr),\epsilon\in V,\mu\in M\bigr\}$ acts on $\mathbb C^2$ properly discontinuous and without fixed points as follows:
\begin{equation}\label{action1}
\big(\begin{smallmatrix}\epsilon &\mu\\ 0 & 1\end{smallmatrix}\big)\bullet(z_1,z_2):=(\epsilon z_1+\mu,\epsilon' z_2+\mu'),
\end{equation}
where $M\cong \mathbb Z\times\mathbb Z\omega_0$ is a rank 2 free module over 
$\mathbb Z$ generated by 1 and another real quadratic irrational number $\omega_0$  and $V\cong \mathbb Z$ generated by a quadratic irrational number $\epsilon>0$. Here $\epsilon'$ (resp. $\mu'$) is the Galois conjugation over $\mathbb Q$ of $\epsilon$ (resp. $\mu'$) and $\epsilon$ also satisfies $\epsilon'=\frac{1}{\epsilon}$. The action of $G(M,V)$ can be restricted onto $\mathbb H^2$ where $\mathbb H$ is the upper half plane.

Let $x_1,y_1$ and $x_2,y_2$ be the real coordinates of $z_1,z_2$ separately and define $$L_c:=\{y_1y_2=c\}\subset \mathbb H\times\mathbb H,$$ where $c$ is a positive constant.  By $\epsilon\epsilon'=1$,  we know that $L_c$ is invariant under the action of $\Gamma$. Now show that the  set $L_c/\Gamma$ is a torus bundle over a circle. Let $x_1,x_2,y_1$ be the coordinates  on set $L_c$. By \eqref{action1}, the action of $\bigl( \begin{smallmatrix}\epsilon^n & \mu\\ 0 & 1\end{smallmatrix}\bigr)$ when restricted on $L_c$ is given by 
$$(x_1,x_2,y_1)\longrightarrow (\epsilon^n x_1+\mu,\epsilon^{-n}x_2+\mu',\epsilon^n y_1)$$
The action of $M$ via the embedding
\begin{equation}\mu\mapsto\Big(\begin{matrix}1&\mu\\ 0 & 1\end{matrix}\Big)\end{equation} on $\mathbb R^2(x_1,x_2)$ realizes $M$ as a lattice of $\mathbb R^2$, so that  the projection of $L_c/\Gamma$ onto $y_1$ identifies $L_c/\Gamma$ with a $\mathbb T^2$ bundle over $\mathbb S_1$.

%By performing the action of $\epsilon$, we can assume $y_1$ lies in $[1,\epsilon)$. Further performing the action of $\mu$, we can assume $(x_1,x_2)$ lie in a fundamental domain of the lattice action. These observation shows that that $(y_1y_2=c)=\mathbb R^3/\Gamma$ is $\mathbb T^2$ bundle on $\mathbb S^1$.

The metric tensor is given by: 
\begin{align*}
g&=\frac{dx_1\otimes dx_1}{y_1^2}+\frac{dy_1\otimes dy_1}{y_1^2}+\frac{dx_2\otimes dx_2}{y_2^2}+\frac{dy_2\otimes dy_2}{y_2^2}\\
&=\frac{dx_1\otimes dx_1}{y_1^2}+\frac{y_1^2dx_2\otimes dx_2}{c^2}+2\frac{dy_1\otimes dy_1}{y_1^2}+\frac{dc\otimes dc}{c^2}-2\frac{dc\otimes dy_1}{cy_1},
\end{align*}
where $c=y_1y_2$. Letting $\hat c=\log c,$ we have 
$$g=\frac{dx_1\otimes dx_1}{y_1^2}+\frac{y_1^2dx_2\otimes dx_2}{e^{2\hat c}}+2\frac{dy_1\otimes dy_1}{y_1^2}+{d\hat c\otimes d\hat c}-2\frac{d\hat c\otimes dy_1}{y_1}$$
Let $\mathbb T^2(y_1)$ be the family of torus defined by the lattice action on $\mathbb R^2 (x_1,x_2),$ which depends on $y_1$. We claim that the diameter of $(\mathbb T^2,g_{\hat c})\rightarrow 0$ when $\hat c\rightarrow \infty,$ where $g_{\hat c}=\frac{dx_1\otimes dx_1}{y_1^2}+\frac{y_1^2dx_2\otimes dx_2}{e^{2\hat c}}.$ Since $y_1$ is bounded, without loss of generality, we may assume $g_{\hat c}=dx_1\otimes dx_1+\frac{dx_2\otimes dx_2}{e^{2\hat c}}.$ We prove the claim by contradiction, if not, we can assume that diam$(\mathbb T^2,g_{\hat c})$ decrease to some constant $\eta\neq 0$. Hence for each $i\in \mathbb Z$, we have point $q_i\in \mathbb T^2$ such that $dist_{g_i}(o,q_i)\ge\frac{\eta}{2}$ We can take an accumulation point $q$ of $q_i$ by the compactness of $\mathbb T^2$. By triangle inequality, we have $dist_{g_i}(o,q)\ge\frac{\eta}{4}.$ This contradicts the fact that $\mathbb Z+\mathbb Z\omega$ is dense in $\mathbb R$, which enables us to show that 
$dist_{g_i}(o,q)\rightarrow 0.$ So if we take a sequence of base point $p_i=({x_1}_i,{x_2}_i,{y_1}_i,c_i)$ in $\cU\setminus p$ with $c_i\rightarrow 0$,  the Gromov-Hausdorff limit of $(\cU,p_j,g_{KE})$ is a flat cylinder $Y$. More precisely,
$Y$ is defined as $\mathbb R^2/\mathbb Z$,
where $$n\bullet(\hat {y_1},\hat c)=(\hat{y_1}+n\log\epsilon, \hat c)$$
and the metric is $$g_Y=2d\hat{y_1}\otimes d\hat{y_1}+d\hat c\otimes d\hat c-d\hat y_1\otimes d\hat c-d\hat c\otimes\hat {y_1}$$ where $\hat y=\log y$.

\end{enumerate}
\end{proof}
\begin{proof}[Proof of theorem \ref{conj} (Surface case)] By Theorem \ref{Dirichelt-bounededgeometry}, when $\omega=\ddbar\rho$ is K\"ahler-Einstein metric with bounded geometry property,  the following Dirichlet problem on $\cU\setminus p$ admits a bounded solution $\varphi$.
\begin{equation}
\begin{cases}&(\omega + \ddbar \varphi ) ^n =e^{\varphi} {\omega^n} \\
&  \varphi|_{\partial\cU}=\psi,
\end{cases} 
\end{equation}
Now we rewrite the equation
$$( \omega + \ddbar \varphi ) ^n =e^{\varphi} {\omega^n}$$ as $$(\ddbar \tilde\varphi ) ^n =e^{\tilde\varphi-\rho} \frac{{\omega^n}}{(\sqrt{-1})^n\nu \wedge \overline \nu}(\sqrt{-1})^n\nu \wedge \overline \nu,$$
where $\tilde\varphi=\rho+\varphi$ and $\nu$ is a local  holomorphic volume form (possibly multi valued). Note that $\log\frac{{\omega^n}}{(\sqrt{-1})^n\nu \wedge \overline \nu}-\rho$ is pluriharmonic function on $\cU\setminus p$.  \\

\noindent\textbf{Claim:} $\log\frac{{\omega^n}}{(\sqrt{-1})^n\nu \wedge \overline \nu}-\rho$ is a sublog function, i.e $\log\frac{{\omega^n}}{(\sqrt{-1})^n\nu \wedge \overline \nu}-\rho\geq \epsilon\log|\sigma_D|_{h_D}+C_\epsilon$. \\

Assuming this claim, we can complete the proof as follows. By the uniqueness Theorem \ref{uniq}, we conclude that $\tilde\varphi+\log\frac{{\omega^n}}{(\sqrt{-1})^n\nu \wedge \overline \nu}-\rho $ coincides with the solution from Theorem \ref{main1} (when the boundary value is matched). Recall that the Riemannian metric induced by $\omega+\ddbar\varphi$ in Theorem \ref{Dirichelt-bounededgeometry} is complete, in particular, this implies that the K\"ahler-Einstein metric associated with the solution from Theorem \ref{main1} is complete towards $p$.

\begin{proof}[Proof of the Claim] We will use the explicit expression of $\rho,\omega, $ and $\nu$ from Lemma \ref{quasi1} to prove the claim case by case. Now we outline the proof. There are essentially 4 types of isolated log canonical (non Klt) singularities in complex dimension 2.
\begin{enumerate}
    \item [(1.1)] When the  exceptional divisor of the minimal resolution  of singularity is an elliptic curve $D$. Then the model metric $\omega=-\ddbar\rho, \rho=\log(-\log|\sigma_D|_{h_D})$. And when pulled back to the resolution, $\nu\wedge\bar\nu$ has a pole $\frac{1}{|\sigma_D|^2_{h_D}}$ along the exceptional divisor. Straightforward calculation shows that  both $\log\frac{\omega^2}{\nu\wedge\bar\nu}$ and $\rho$ are  sub log functions. 
    \item [(1.2)]The singularity is finite quotient of case (1.1) above. When pulled back to the finite cover, $\log\frac{{\omega^n}}{(\sqrt{-1})^n\nu \wedge \overline \nu}-\rho$,  by proof of case 1 above, is a bounded function. Here we have used the fact that $\omega$ is invariant under the finite quotient, and when pull backed to finite cover, $\nu$ is still a generator of holomorphic volume form near the singularity point $p$. 
   \item [(2.1)]When the  exceptional divisors of the minimal resolution  of singularity is circle of rational curves $D_i, 0\le i\le m$.  Now let us focus on the metric behaviour of $\omega$ near the intersection of two exceptional divisors. Without loss of generality, assume locally $D_k=\{u_k=0\},D_{k+1}=\{v_k=0\}$, then the metric $\omega=\ddbar\rho,\rho=-\log(\log|u_k^\alpha v_k^\beta|\log|u_k^{\alpha'} v_k^{\beta'}|)$
   
   $$\omega=\frac{\sqrt{-1}\partial\log|u_k^\alpha v_k^\beta|\wedge\bar\partial\log|u_k^\alpha v_k^\beta|}{(\log|u_k^\alpha v_k^\beta|)^2}+\frac{\sqrt{-1}\partial\log|u_k^{\alpha'} v_k^{\beta'}|\wedge\bar\partial\log|u_k^{\alpha'} v_k^{\beta'}|}{(\log|u_k^{\alpha'} v_k^{\beta'}|)^2}$$
   And when pulled back to the resolution, $\nu\wedge\bar\nu$ has a pole $\frac{1}{|u_kv_k|^2}$ along the exceptional divisor. Then direct calculation shows that  both $\log\frac{\omega^2}{(\sqrt{-1})^n\nu\wedge\bar\nu}$ and $\rho$ are sublog functions.
   \item [(2.2)]When the singularity is finite quotient of case (2.1), then we can argue as case (1.2) to complete the proof.
\end{enumerate}
\end{proof}

We end this section by providing an example of a family of canonical polarized surface with the central fiber equipped with log canonical singularity satisfying bounded geometry property.
\begin{example}\label{ex} Degeneration of Godeaux surfaces
A surface $X$ is called a Godeaux surface if $\pi_1(X)=\mathbb Z_5$ and universal cover is quintic hypersurface. A explicit construction could be as follows: Define $\mathbb Z_5$ on $\mathbb P^3$ in the following way:
$$\rho\bullet(X_0,X_1,X_2,X_3)=(X_0,\rho X_1,\rho^2X_2.\rho^3 X_3)$$
Then there exists quintics (in $\mathbb P^3$) invariant and fixed point free under the $\mathbb Z_5$ action with 5 non degenerate triple points and no other singularities by a dimension count argument. (cf. \cite[pg.\ 135]{Ulf}). Then the $\mathbb Z_5$ quotient will give a family of Godeaux with central fiber a canonical polarized variety coupled with a single simple elliptic singularity (cone over elliptic curve). 
\end{example}
\end{proof}

\section{Discussion of conjecture \ref{conj1} for Smoothable isolated log canonical singularities}
In this section, we give a discussion of conjecture \ref{conj1} when the isolated log canonical singularity is smoothable.
First we set up the question. 

Let  $\pi: \mathcal{X} \rightarrow B$ be a smoothing of an isolated log canonical (non-log terminal) singularity $(\mathcal X_0,p)$ over $B\in \mathbb{C}$ such that $\mathcal X_t= \pi^{-1}(t) $ is smooth for  $t\neq 0$. We  assume that the total space $\mathcal{X}$ has at worst canonical singularities and the canonical divisor $K_{\mathcal{X}}$ is ample. For simplicity, we assume $K_{\mathcal{X}/B}$ is Cartier and let  $\nu$ be  a local generator $\nu$ of the bundle $K_{\mathcal{X}/B}$ such that  $\nu_t = \nu |_{\mathcal{X}_t}$ is a local holomorphic $n$ form on $\mathcal{X}_t$, where $\dim \mathcal{X}_t = n$ for $t\in B$. $\mathcal X$ can locally be embedded in $
\mathbb C^{N+1}$ and we apply a pluriharmonic function $\rho$ on $\mathbb C^{N+1}$ to  define $\cU:=\{\rho\le 0 \}\cap \mathcal{X}$ so that $p\in \cU$.  After perturbation, we can always assume $\partial \cU$ is smooth and $d\rho \neq 0$  on $\partial \cU$. By continuity, we may assume $\partial \cU_t$ is smooth for $t$ sufficiently close to $0$, where $\cU_t = \cU\cap \mathcal{X}_t$ and $\rho_t = \rho|_{\mathcal{X}_t}$. 

We let $\psi$ be a smooth function in an open neighborhood of  $\partial \cU$ and let $\psi_t$ be its restriction to $\partial \cU_t$. 
%$$\chi =\frac{1}{m}  \ddbar \log \left( \sum_{j=0}^{N} |\eta_j|^2 \right).$$ 
%
%$$\omega=\ddbar\rho, ~\omega_t = \omega|_{\cU_t}=\ddbar\rho_t $$ be a family smooth K\"ahler form on $\mathcal{X}_t$  for $t\in B$. For example, we can choose $\rho= \sum_{j=1}^{N+1} |z_j|^2 - \epsilon$ for sufficiently small $\epsilon>0$ and so $\omega$ is simply the standard Euclidean metric on $\mathbb{C}^{N+1}$. 
 We also let  $\Omega$ be the real-valued smooth $(n,n)$-volume measure  on $\mathcal{X}$ defined by 
$$\Omega =  (\sqrt{-1})^n \nu\wedge \bar \nu $$ and let 
\begin{equation}\label{V}\Omega_t = \Omega|_{\mathcal{X}_t}= (\sqrt{-1})^n \nu_t\wedge \bar \nu_t\end{equation} be the restriction of $\Omega$ on $\cU_t$. Immediately we have  
$$\ddbar \log \Omega = 0, ~~\ddbar\log \Omega_t = 0$$
for $t\in B$.

 By Theorem \ref{main1},  there is a unique solution $\varphi_t$ to the following Dirichlet problem for Monge-Amp\`ere equation on $\cU_t$ 
\begin{equation}\label{familyscmp}
\begin{cases}
&(\ddbar \varphi_t)^n = e^{\varphi_t}\Omega_t  \\
&\varphi_t|_{\partial\cU_t}=\psi_t
\end{cases}
\end{equation}
for each $t\in B$. If we let $\omega_t = \ddbar \varphi_t$, then the corresponding Kahler metric $g_t$ satisfies the Kahler-Einstein equation 
$$Ric(g_t) = - g_t$$ on $\cU_t$ for $t\in B^*$ and $\cU_0\setminus \{p\}$ for $t=0$. 

In \cite{S4}, the third author combines semi-stable reduction and maximum principle to derive uniform estimates with barrier for $\varphi_t$ for $t\in B$ when $\mathcal{X}$ is a stable family of canonical polarized variety. And then apply the local version of partial $C^0$  estimate (cf. \cite{DS1,DS2}) initiated by \cite{Tian} for noncollapsed polarized K\"ahler-Einstein manifold and some other tools to show that the K\"ahler-Einstein metric $g_0$ on $\mathcal{X}_0$ is complete towards the local canonical locus.
   
In our local situation, now we have a degeneration of polarized K\"ahler-Einstein manifold $$(\cU_t,g_t,L_t,h_t)$$ with boundary, where $L_t$ is a trivial line bundle and  $h_t:=e^{-\varphi_t}$ is the hermitian metric on $L_t$ whose curvature is the K\"ahler form associated to the K\"ahler-Einstein metric  $g_t$ . 
%is   following data $h_t = ((\omega^{KE}_t)^n)^{-1}= (e^{\varphi_t}\Omega_t )^{-1}$ be the hermitian metric on $\mathcal{X}_t$ for $t\in B^*$, where $\omega^{KE}_t$ is K\"ahler-Einstein form associated to the Kahelr-Einstein metric $g_t$ on $\mathcal{X}_t$, $\Omega_t$ and $\varphi_t$ are defined by equation \ref{familyscmp}. 
Notice that $\cU_t, t\in B^*$ is a strongly pseudoconvex domain, so we can solve  $\bar\partial$ equation with $L^2$ estimate on $\Omega_t$. So we believe that the techniques developed in \cite{S4} and the local version of the partial $C^0$-estimate (cf. \cite[Theorem 1.1]{DS2}) can still be applied.

\bigskip

\begin{appendices}
\section{Appendix}
\begin{lemma}\label{comparisonp} Let $(\cU,\omega)$ be a K\"ahler manifold with boundary (the boundary can possibly have several components). If $\varphi,\psi$ are two continuous $\omega-$ PSH function on $\bar\cU$ satisfying $\varphi>\psi$ on $\partial\cU$, then for $\Omega:=\{\varphi<\psi\}$ we have 
$$\int_\Omega\omega^n_\psi\leq\int_\Omega\omega^n_\varphi.$$
\end{lemma}
\begin{proof}
Suppose first that $\varphi,\psi$ and the boundary $\Omega$ are smooth. Set $\varphi_t=\max(\varphi+t,\psi), t>0$. Then close to $\partial\Omega$, we have $\varphi_t=\varphi+t$. Define the closed current
$$T_t=\sum_{k=1}^{k=n}{{n}\choose {k}}(\ddbar\varphi_t)^{k-1}\wedge\omega^{n-k}$$
and set $T=\lim_{t\to 0}T_t$. By Stokes's theorem
\begin{align*}
    \int_{\Omega}\omega_{\varphi_t}^n&=\int_{\Omega}\ddbar\varphi_t\wedge T_t+\omega^n=\int_{\partial\Omega}d^c\varphi_t\wedge T_t+\int_\Omega\omega^n\\
    &=\int_{\partial\Omega}d^c\varphi\wedge T+\int_{\Omega}\omega^n=\int_\Omega\omega^n_{\varphi}
\end{align*}
Since $\varphi_t\to \psi$ in $\Omega$ as $t\to 0$, we get by applying the convergence theorem that $\omega^n_{\varphi_t}\to\omega^n_{\psi}$ in $\Omega$. Hence for a test function $\chi$ in $\Omega$ with $0\leq\chi\leq 1$ we get
$$\int_{\Omega}\chi\omega^n_{\psi}=\lim_{t\to 0}\int_\Omega\chi\omega^n_{\varphi_t}\leq\liminf_{t\to 0}\int_\Omega\omega^n_{\varphi_t}.$$
So
$$\int_\Omega\omega^n_\psi\leq\liminf_{t\to 0}\int_\Omega\omega^n_{\varphi_t}=\int_\Omega\omega^n_\varphi,$$ which completes the proof for smooth functions. Now we suppse that $\varphi,\psi$ satisfy the extra assumption 
\begin{equation}\label{a}
    \ddbar\varphi\geq(\delta-1)\omega,\ddbar\psi\geq(\delta-1)\omega,
\end{equation}
for some $\delta>0$. Then by \cite[Theorem 2]{BlockiK}, we can find two sequence of $\omega-$ PSH function $\varphi_j$ and $\psi_j$ on $\bar\Omega$ converging uniformly to $\varphi,\psi$ respectively. Given a compact set $K\subset\Omega$, we find $t>0$ and a positive integer $j_0$ such that 
$$K\subset\Omega(t,j):=\{\varphi_j<\psi_j-t\}\subset\Omega$$
for $j>j_0$ and the boundary of $\Omega(t,j)$ is smooth (using Sard's theorem). Now we have 
$$\int_K\omega_\psi^n\leq\liminf_{j\to\infty}\int_{\Omega(t,j)}\omega^n_{\psi_j}\leq\liminf_{j\to\infty}\int_{\Omega(t,j)}\omega^n_{\varphi_j}\leq\int_\Omega\omega^n_{\varphi},$$
where the second inequality is due to the first part of the proof. We still need to get rid of the assumption of (\ref{a}). Note that for fixed $t\in(0,1)$ and $\omega-$ PSH functions $\varphi,\psi$, the functions $t\varphi,t\psi$ satisfy (\ref{a}) for some $\delta>0$. For a fixed compact set $K\subset\Omega$ and constant $t\in (0,1)$, we may choose $\delta>0$ sufficiently small such that $K\subset\Omega(\delta,t):=\{\varphi<\psi-\frac{\delta}{t}\}$. Now we have 
$$\int_K\omega_\psi^n\leq\liminf_{t\to 1}\int_{\Omega(t,\delta)}\omega^n_{t\psi}\leq\liminf_{t\to 1}\int_{\Omega(t,\delta)}\omega^n_{t\varphi}\leq\int_\Omega\omega^n_{\varphi},$$
To complete the proof it is enough to consider an exhaustion sequence of compact subsets of $\Omega$.
\end{proof}

\end{appendices}

\end{document}